\documentclass[11pt]{amsart}
\usepackage{amsmath, amsthm, amscd, amsfonts, amssymb, graphicx, color,float,pgf,tikz}
\usepackage{amssymb,fontenc}
\usepackage{latexsym,wasysym,mathrsfs}
\usepackage{hyperref}
\usetikzlibrary{arrows}
\usepackage[square,numbers,sort&compress]{natbib}
\makeatletter \oddsidemargin.9375in \evensidemargin \oddsidemargin
\newtheorem{theorem}{Theorem}[section]

\theoremstyle{definition}
\newtheorem{definition}[theorem]{Definition}

\newtheorem{remark}[theorem]{Remark}
\numberwithin{equation}{section}

\newcommand{\be}{\begin{equation}}
\newcommand{\ee}{\end{equation}}

\numberwithin{equation}{section}

\usepackage{etoolbox}
\makeatletter
\patchcmd{\@settitle}{\uppercasenonmath\@title}{}{}{}
\patchcmd{\@setauthors}{\MakeUppercase}{}{}{}
\makeatother

\makeatletter
\@namedef{subjclassname@2020}{%
  \textup{2020} Mathematics Subject Classification}
\makeatother

\allowdisplaybreaks
\begin{document}

\noindent {\tiny Sahand Communications in Mathematical Analysis (SCMA) Vol. $\cdots$ No. $\cdots$(2021), $\cdots$-$\cdots$\\
\url{http://scma.maragheh.ac.ir}}\\
DOI: 10.22130/scma.2021.524252.909   \\[0.50in]


\title[Introduction of frame in tensor product of \,$n$-Hilbert spaces]{Introduction of frame in tensor product of \,$n$-Hilbert spaces}

\author[P. Ghosh and T. K. Samanta]{Prasenjit Ghosh$^1$$^{*}$ and Tapas Kumar Samanta$^2$}

\address{ $^{1}$ Department of Pure Mathematics, University of Calcutta, 35, Ballygunge Circular Road, Kolkata, 700019, West Bengal, India.}
\email{prasenjitpuremath@gmail.com}

\address{ $^{2}$ Department of Mathematics, Uluberia College, Uluberia, Howrah, 711315,  West Bengal, India.}
\email{mumpu$_{-}$tapas5@yahoo.co.in}



\keywords{Frame, Dual frame, Tensor product of Hilbert spaces, $n$-normed space, $n$-Hilbert space.\\
\indent Received: dd mmmm yyyy,    Accepted: dd mmmm yyyy.
\\
\indent $^{*}$ Corresponding author}
\maketitle
\hrule width \hsize \kern 1mm


\begin{abstract}
We study the concept of frame in tensor product of \,$n$-Hilbert spaces as tensor product of \,$n$-Hilbert spaces is again a \,$n$-Hilbert space.\,We generalize some of the known results about bases to frames in this new Hilbert space.\,A relationship between frame and bounded linear operator in tensor product of \,$n$-Hilbert spaces is studied.\,Finally,\;the dual frame in tensor product of \,$n$-Hilbert spaces is discussed.
\end{abstract}
\maketitle
\vspace{0.1in}
\hrule width \hsize \kern 1mm

\section{Introduction}
There are several techniques to rebuild signals by using a family of elementary signals.\,One of these techniques was given by D. Gabor in 1946 \cite{Gabor}.\,Thereafter frame in Hilbert space has been developed by Duffin and Schaeffer \cite{Duffin}.\,It is very useful to study nonharmonic Fourier series, i.e., sequences of the type \,$\left\{\,e^{\,i\,\lambda_{n}\,x}\,\right\}_{n \,\in\, \mathbb{Z}}$, where \,$\left\{\,\lambda_{n}\,\right\}_{n \,\in\, \mathbb{Z}}$\, is a family of real or complex numbers.\;Daubechies et al, connected frames with Gabor system and Wavelet in 1986  \cite{Daubechies}.\,To translate informations, frames are more flexible tools than bases and which makes them very useful in signal processing \cite{F}, coding and communications \cite{H}, system modeling \cite{E}, filter bank theory \cite{BH} etc.\,In recent times, many generalizations of frames have been appeared.\,Some of them are \,$g$-frame \cite{Sun}, fusion frame \cite{Kutyniok} and \,$g$-fusion frame \cite{Ahmadi} etc.\;P.\,Ghosh and T.\,K.\,Samanta studied the stability of dual \,$g$-fusion frames and generalized atomic systems for operators in Hilbert spaces \cite{P, Ghosh}.

In 1970, Diminnie et al, introduced the concept of \,$2$-inner product space \cite{Diminnie}.\;A generalization of \,$2$-inner product space for \,$n \,\geq\, 2$\, was developed by A.\,Misiak in 1989 \cite{Misiak}.\;The basic concepts of tensor product of Hilbert spaces were presented by S.\,Rabinson \cite{S}.

In this paper, we study frame in the tensor product of \,$n$-Hilbert spaces and establish some of its properties.\;We note the result that in tensor product of \,$n$-Hilbert spaces, an image of a frame under a bounded linear operator is a frame if and only if the operator is invertible.\;Finally, dual of a frame in tensor product of \,$n$-Hilbert spaces is described. 

Throughout this paper,\,$X$\, denotes separable Hilbert space associated with the inner product \,$\left<\,\cdot \,,\, \cdot\,\right>_{1}$\, and \,$l^{\,2}(\,\mathbb{N}\,)$\; denote the space of square summable scalar-valued sequences with index set of natural numbers \,$\mathbb{N}$.

\section{Preliminaries}

\begin{definition}\cite{Christensen}
A sequence \,$\left\{\,p_{\,i}\,\right\}_{i \,=\, 1}^{\infty} \,\subseteq\, X$\, is said to be a frame for \,$X$\, if there exist positive constants \,$A,\, B$\, such that
\[ A\; \|\,p\,\|_{1}^{\,2} \,\leq\, \;\sum\limits_{i \,=\, 1}^{\infty}\;  \left|\ \left <\,p,\, p_{\,i} \, \right >_{1}\,\right|^{\,2} \,\leq\, B \;\|\,p\,\|_{1}^{\,2};\; \;\forall\; p \,\in\, X.\]
The constants \,$A$\, and \,$B$\, are called frame bounds.
\end{definition}

\begin{definition}\cite{Christensen}
Let \,$\{\,p_{\,i}\,\}_{i \,=\, 1}^{\infty}$\, be a frame for \,$X$.\;The synthesis operator, \,$ T \,:\, l^{\,2}\,(\,\mathbb{N}\,) \,\to\, X$, defined by \,$T\,\{\,c_{i}\,\} \,=\, \sum\limits_{i \,=\, 1}^{\infty} \;c_{\,i}\,p_{\,i}$\, is bounded operator and its adjoint, called the analysis operator, is given by \,$T^{\,\ast}\,p \,=\, \left \{\, \left <\,p,\, p_{i} \,\right >_{1}\,\right \}_{i \,=\, 1}^{\infty}$.\;The frame operator \,$S \,:\, X \,\to\, X$, is given by \\$S\,p \,=\, T\,T^{\,\ast}\,p \,=\, \sum\limits^{\infty}_{i \,=\, 1}\; \left <\,p,\, p_{\,i} \, \right >_{1}\,p_{\,i}$, for all \,$p \,\in\, X$.
\end{definition}

\begin{definition}\cite{Christensen}
A frame \,$\left\{\,q_{\,i}\,\right\}_{i \,=\, 1}^{\infty}$\, is said to be a dual of a frame \,$\left\{\,p_{\,i}\,\right\}_{i \,=\, 1}^{\infty}$\, for \,$X$\, if \,$p \,=\, \sum\limits^{\infty}_{i \,=\, 1}\; \left <\,p,\, q_{\,i} \, \right >_{1}\,p_{\,i}$, for all \,$p \,\in\, X$. 
\end{definition}

The tensor product of Hilbert spaces are introduced by several ways and it is a certain linear space of operators which was represented by Folland in \cite{Folland}, Kadison and Ringrose in \cite{Kadison}.

\begin{definition}\cite{Upender}\label{def0.001}
Let \,$\left(\,Y,\, \left<\,\cdot,\,\cdot\,\right>_{2}\,\right)$\, be a Hilbert space.\,Then the tensor product of Hilbert spaces \,$X$\, and \,$Y$\, is denoted by \,$X \,\otimes\, Y$\, and it is defined to be an inner product space with respect to the inner product: 
\[\left<\,p \,\otimes\, q,\, p^{\,\prime} \,\otimes\, q^{\,\prime}\,\right> \,=\, \left<\,p,\, p^{\,\prime}\,\right>_{1}\;\left<\,q,\, q^{\,\prime}\,\right>_{2}, \;\text{for all}\; p,\, p^{\,\prime} \,\in\, X; \;q,\, q^{\,\prime} \,\in\, Y.\]\,The norm on \,$X \,\otimes\, Y$\, is given by 
\[\left\|\,p \,\otimes\, q\,\right\| \,=\, \|\,p\,\|_{\,1}\;\|\,q\,\|_{\,2}, \;\text{for all}\; p \,\in\, X\; \;\text{and}\; \,q \,\in\, Y.\]\,The space \,$X \,\otimes\, Y$\, is complete with respect to the above inner product.\;Therefore the space \,$X \,\otimes\, Y$\, is a Hilbert space.     
\end{definition} 

Tensor product of operators \,$Q \,\in\, \mathcal{B}\,(\,X\,)$\, and \,$T \,\in\, \mathcal{B}\,(\,Y\,)$, is denoted by \,$Q \,\otimes\, T$\, and defined as \,$\left(\,Q \,\otimes\, T\,\right)\,A \,=\, Q\,A\,T^{\,\ast}, \;\text{for all}\; A \,\in\, X \,\otimes\, Y$.
It can be easily verified that \,$Q \,\otimes\, T \,\in\, \mathcal{B}\,(\,X \,\otimes\, Y\,)$\, \cite{Folland}.

\begin{theorem}\cite{Folland}\label{th1.1}
Suppose \,$Q,\, Q^{\prime} \,\in\, \mathcal{B}\,(\,X\,)$\, and \,$T,\, T^{\prime} \,\in\, \mathcal{B}\,(\,Y\,)$.\,Then \begin{itemize}
\item[$(i)$]\,$Q \,\otimes\, T \,\in\, \mathcal{B}\,(\,X \,\otimes\, Y\,)$\, and \,$\left\|\,Q \,\otimes\, T\,\right\| \,=\, \|\,Q\,\|\; \|\,T\,\|$.
\item[$(ii)$] \,$\left(\,Q \,\otimes\, T\,\right)\,(\,f \,\otimes\, g\,) \,=\, Q\,f \,\otimes\, T\,g$\, for all \,$f \,\in\, X,\, g \,\in\, Y$.
\item[$(iii)$] $\left(\,Q \,\otimes\, T\,\right)\,\left(\,Q^{\,\prime} \,\otimes\, T^{\,\prime}\,\right) \,=\, (\,Q\,Q^{\,\prime}\,) \,\otimes\, (\,T\,T^{\,\prime}\,)$. 
\item[$(iv)$] \,$Q \,\otimes\, T$\, is invertible if and only if \,$Q$\, and \,$T$\, are invertible, in which case \,$\left(\,Q \,\otimes\, T\,\right)^{\,-\, 1} \,=\, \left(\,Q^{\,-\, 1} \,\otimes\, T^{\,-\, 1}\,\right)$.
\item[$(v)$] \,$\left(\,Q \,\otimes\, T\,\right)^{\,\ast} \,=\, \left(\,Q^{\,\ast} \,\otimes\, T^{\,\ast}\,\right)$.    
\end{itemize}
\end{theorem}

\begin{definition}\cite{Mashadi}
A real valued function \,$\left\|\,\cdot,\, \cdots,\, \cdot \,\right\| \,:\, H^{\,n} \,\to\, \mathbb{R}$\, satisfying the following properties:
\begin{itemize}
\item[$(i)$]\;\; $\left\|\,x_{\,1},\, x_{\,2},\, \cdots,\, x_{\,n}\,\right\| \,=\, 0$\; if and only if \,$x_{\,1},\, \cdots,\, x_{\,n}$\; are linearly dependent,
\item[$(ii)$]\;\;\; $\left\|\,x_{\,1},\, x_{\,2},\, \cdots,\, x_{\,n}\,\right\|$\, is invariant under permutations of \,$x_{1},\, \cdots, x_{n}$,
\item[$(iii)$]\;\;\; $\left\|\,\alpha\,x_{\,1},\, x_{\,2},\, \cdots,\, x_{\,n}\,\right\| \,=\, |\,\alpha\,|\,\left\|\,x_{\,1},\, x_{\,2},\, \cdots,\, x_{\,n}\,\right\|$, \,$\alpha \,\in\, \mathbb{K}$,
\item[$(iv)$]\;\; $\left\|\,x \,+\, y,\, x_{\,2},\, \cdots,\, x_{\,n}\,\right\| \,\leq\, \left\|\,x,\, x_{\,2},\, \cdots,\, x_{\,n}\,\right\| \,+\,  \left\|\,y,\, x_{\,2},\, \cdots,\, x_{\,n}\,\right\|$,
\end{itemize}
for all \,$x_{\,1},\, x_{\,2},\, \cdots,\, x_{\,n},\,x,\, y \,\in\, H$, is called \,$n$-norm on \,$H$.\,A linear space \,$H$, together with a \,$n$-norm \,$\left\|\,\cdot,\, \cdots,\, \cdot \,\right\|$, is called a linear\;$n$-normed space. 
\end{definition}

\begin{definition}\cite{Misiak}\label{defnn1}
Let \,$n \,\in\, \mathbb{N}$\; and \,$H$\, be a linear space of dimension greater than or equal to \,$n$\, over the field \,$\mathbb{K}$, where \,$\mathbb{K}$\, is the real or complex numbers field.\;A function \,$\left<\,\cdot,\, \cdot \,|\, \cdot,\, \cdots,\, \cdot\,\right> \,:\, H^{n \,+\, 1} \,\to\, \mathbb{K}$\, satisfying the following five properties: 
\[(i)\left<\,x_{\,1},\, x_{\,1} \,|\, x_{\,2},\, \cdots,\, x_{\,n} \,\right> \,\geq\,  0\; \text{and} \;\left<\,x_{\,1},\, x_{\,1} \,|\, x_{\,2},\, \cdots,\, x_{\,n} \,\right> \,=\,  0\hspace{1cm}\]
 if and only if \,$x_{\,1},\, x_{\,2},\, \cdots,\, x_{\,n}$\, are linearly dependent,
\[(ii)\left<\,x,\, y\,|\,x_{\,2},\, \cdots,\, x_{\,n} \,\right> \,=\, \left<\,x,\, y\,|\,x_{\,i_{\,2}},\, \cdots,\, x_{\,i_{\,n}} \,\right>\text{ for every permutations}\]
$\left(\, i_{\,2},\, \cdots,\, i_{\,n} \,\right)$\, of \,$\left(\, 2,\, \cdots,\, n \,\right)$,
\[(iii)\left<\,x,\, y \,|\, x_{\,2},\, \cdots,\, x_{\,n} \,\right> \,=\, \overline{\left<\,y,\, x \,|\, x_{\,2},\, \cdots,\, x_{\,n} \,\right> },\hspace{4cm}\]
\[(iv)\left<\,\alpha\,x,\, y \,|\, x_{\,2},\, \cdots,\, x_{\,n} \,\right> \,=\, \alpha \,\left<\,x,\, y \,|\, x_{\,2},\, \cdots,\, x_{\,n}\,\right>, \,\alpha \,\in\, \mathbb{K}\hspace{4cm}\]
\[(v)\left<\,x \,+\, y,\, z\,|\,x_{2},\, \cdots,\, x_{n} \,\right> \,=\, \left<\,x,\, z\,|\,x_{2},\, \cdots,\, x_{n} \,\right> \,+\,  \left<\,y,\, z\,|\,x_{2},\, \cdots,\, x_{n} \,\right>,\]
for all \,$x,\, y,\, x_{\,1},\, x_{\,2},\, \cdots,\, x_{\,n} \,\in\, H$, is called an \,$n$-inner product on \,$X$, and the pair \,$\left(\,X,\, \left<\,\cdot,\, \cdot \,|\, \cdot,\, \cdots,\, \cdot\,\right>\,\right)$\, is called an \,$n$-inner product space.
\end{definition}

\begin{theorem}\cite{Misiak}
Let \,$H$\, be an \,$n$-inner product space.\,Then 
\[\left|\,\left<\,x,\, y \,|\, x_{\,2},\,  \cdots,\, x_{\,n}\,\right>\,\right| \,\leq\, \left\|\,x,\, x_{\,2},\, \cdots,\, x_{\,n}\,\right\|\, \left\|\,y,\, x_{\,2},\, \cdots,\, x_{\,n}\,\right\|,\]
for all \,$x,\, y,\, x_{\,2},\, \cdots,\, x_{\,n} \,\in\, H$, where  \,\[\left \|\,x_{\,1},\, x_{\,2},\, \cdots,\, x_{\,n}\,\right\| \,=\, \sqrt{\left <\,x_{\,1},\, x_{\,1} \,|\, x_{\,2},\,  \cdots,\, x_{\,n}\,\right>},\] is called Cauchy-Schwarz inequality.
\end{theorem}

\begin{theorem}\cite{Misiak}
Let \,$H$\, be an \,$n$-inner product space.\,Then 
\[\left \|\,x_{\,1},\, x_{\,2},\, \cdots,\, x_{\,n}\,\right\| \,=\, \sqrt{\left <\,x_{\,1},\, x_{\,1} \,|\, x_{\,2},\,  \cdots,\, x_{\,n}\,\right>}\] defines a n-norm for which 
\[\left <\,x,\, y \,|\, x_{\,2},\,  \cdots,\, x_{\,n}\,\right>\hspace{6.5cm}\]
\[ \,=\,\dfrac{\,1}{\,4}\, \left(\,\|\,x \,+\, y,\, x_{\,2},\, \cdots,\, x_{\,n}\,\|^{\,2} \,-\, \|\,x \,-\, y,\, x_{\,2},\, \cdots,\, x_{\,n}\,\|^{\,2}\,\right),\]and 
\[\|\,x \,+\, y,\, x_{\,2},\, \cdots,\, x_{\,n}\,\|^{\,2} \,+\, \|\,x \,-\, y,\, x_{\,2},\, \cdots,\, x_{\,n}\,\|^{\,2}\]
\[ \,=\, 2\, \left(\,\|\,x,\, x_{\,2},\, \cdots,\, x_{\,n}\,\|^{\,2} \,+\, \|\,y,\, x_{\,2},\, \cdots,\, x_{\,n}\,\|^{\,2} \,\right)\] 
hold for all \,$x,\, y,\, x_{\,1},\, x_{\,2},\, \cdots,\, x_{\,n} \,\in\, H$.
\end{theorem}

\begin{definition}\cite{Mashadi}
A sequence \,$\{\,x_{\,k}\,\}$\, in linear\;$n$-normed space \,$H$\, is said to be convergent to \,$x \,\in\, H$\, if 
\[\lim\limits_{k \to \infty}\,\left\|\,x_{\,k} \,-\, x,\, e_{\,2},\, \cdots,\, e_{\,n} \,\right\| \,=\, 0\]
for every \,$ e_{\,2},\, \cdots,\, e_{\,n} \,\in\, H$\, and it is called a Cauchy sequence if 
\[\lim\limits_{l,\, k \,\to\, \infty}\,\left \|\,x_{l} \,-\, x_{\,k},\, e_{\,2},\, \cdots,\, e_{\,n}\,\right\| \,=\, 0\]
for every \,$ e_{\,2},\, \cdots,\, e_{\,n} \,\in\, H$.\;The space \,$H$\, is said to be complete if every Cauchy sequence in this space is convergent in \,$H$.\;A \,$n$-inner product space is called \,$n$-Hilbert space if it is complete with respect to its induce norm.
\end{definition}

\section{Frame in tensor product of $n$-Hilbert spaces}

\begin{definition}\label{def0.1}
Let \,$H$\, be a \,$n$-Hilbert space and \,$a_{\,2},\, \cdots,\, a_{\,n} \,\in\, H$.\;A sequence \,$\left\{\,p_{\,i}\,\right\}^{\,\infty}_{\,i \,=\, 1}$\, in \,$H$\, is said to be a frame associated to \,$\left(\,a_{2},\, \cdots,\, a_{n}\,\right)$\, for \,$H$\, if there exist constants \,$0 \,<\, A \,\leq\, B \,<\, \infty$\,  such that
\begin{equation}\label{eqp1}
A  \left\|\,p,\, a_{\,2},\, \cdots,\, a_{\,n}\,\right\|^{2} \leq \sum\limits^{\infty}_{i \,=\, 1}\left|\left<\,p,\, p_{\,i} \,|\, a_{\,2},\, \cdots,\, a_{\,n}\,\right>\right|^{2} \leq B\, \left\|\,p,\, a_{2},\, \cdots,\, a_{n}\,\right\|^{2}
\end{equation}
for all \,$p \,\in\, H$.\,The constants \,$A$\, and \,$B$\, are called frame bounds.\;If \,$\left\{\,p_{\,i}\,\right\}^{\,\infty}_{\,i \,=\, 1}$\; satisfies the right hand inequality of (\ref{eqp1}), it is called a Bessel sequence associated to \,$\left(\,a_{\,2},\, \cdots,\, a_{\,n}\,\right)$\, in \,$H$\, with bound \,$B$.
\end{definition}
By Cauchy-Schwarz inequality and the condition \,$(i)$\, of the definition \ref{defnn1}, we may assume that every \,$p_{\,i}$\, and \,$a_{\,2},\, \cdots,\, a_{\,n}$\, are linearly independent.

Consider \,$F \,=\, \left\{\,a_{\,2},\, a_{\,3},\, \cdots,\, a_{\,n}\,\right\}$, where \,$a_{\,2},\, a_{\,3},\, \cdots,\,a_{\,n}$\, are fixed elements in \,$H$ and \,$L_{F}$\, denote the linear subspace of \,$H$\, spanned by \,$F$.\;The quotient space \,$H \,/\, L_{F}$\, is a normed linear space with respect to the norm, \,$\left\|\,p \,+\, L_{F}\,\right\|_{F} \,=\, \left\|\,p,\, a_{\,2},\,  \cdots,\, a_{\,n}\,\right\|$\, for every \,$p \,\in\, H$.\;Let \,$M_{F}$\, be the algebraic complement of \,$L_{F}$, then \,$H \,=\, L_{F} \,\oplus\, M_{F}$.\;Define   
\[\left<\,p,\, q\,\right>_{F} \,=\, \left<\,p,\, q \,|\, a_{\,2},\, \cdots,\, a_{\,n}\,\right>\; \;\text{on}\; \;H.\]
Then \,$\left<\,\cdot,\, \cdot\,\right>_{F}$\, is a semi-inner product on \,$H$\, and this induces an inner product on the quotient space \,$H \,/\, L_{F}$\; which is given by
\[\left<\,p \,+\, L_{F},\, q \,+\, L_{F}\,\right>_{F} \,=\, \left<\,p,\, q\,\right>_{F} \,=\, \left<\,p,\, q \,|\, a_{\,2},\,  \cdots,\, a_{\,n} \,\right>;\; \forall\; p,\, q \,\in\, H.\]
By identifying \,$H \,/\, L_{F}$\; with \,$M_{F}$\; in an obvious way, we obtain an inner product on \,$M_{F}$.\;Now, for every \,$p \,\in\, M_{F}$, we define \,$\|\,p\,\|_{F} \,=\, \sqrt{\left<\,p,\, p \,\right>_{F}}$\, and it can be easily verify that \,$\left(\,M_{F},\, \|\,\cdot\,\|_{F}\,\right)$\, is a norm space.\;Consider \,$H_{F}$\, as the completion of the inner product space \,$M_{F}$.

\begin{theorem}\label{th2}
Let \,$H$\, be a n-Hilbert space.\,Then \,$\left\{\,p_{\,i}\,\right\}^{\,\infty}_{\,i \,=\, 1} \,\subseteq\, H$\, is a frame associated to \,$\left(\,a_{\,2},\, \cdots,\, a_{\,n}\,\right)$\, with bounds \,$A,\,B$\, if and only if it is a frame for the Hilbert space \,$H_{F}$\, with bounds \,$A,\,B$.
\end{theorem}

\begin{proof}
Let us consider that \,$\left\{\,p_{\,i}\,\right\}^{\,\infty}_{\,i \,=\, 1} \,\subseteq\, H$\, be a frame associated to \,$\left(\,a_{\,2},\, \cdots,\, a_{\,n}\,\right)$\, for \,$H$\, with bounds \,$A,\,B$.\,Then the inequality (\ref{eqp1}) can be written as 
\[ A\,\|\,p\,\|_{F}^{\,2} \,\leq\, \;\sum\limits_{i \,=\, 1}^{\infty}\;  \left|\ \left <\,p,\, p_{\,i} \, \right >_{F}\,\right|^{\,2} \,\leq\, B \;\|\,p\,\|_{F}^{\,2};\; \;\forall\; p \,\in\, M_{F}.\]
This shows that \,$\left\{\,p_{\,i}\,\right\}^{\,\infty}_{\,i \,=\, 1}$\, is a frame for \,$M_{F}$.\,By the Lemma $5.1.2$\, of \cite{Christensen}, the sequence \,$\left\{\,p_{\,i}\,\right\}^{\,\infty}_{\,i \,=\, 1}$\, is also a frame for \,$H_{F}$\, with the same bounds.

Converse part is obvious.
\end{proof}

\begin{definition}\label{def0.0001}
Let \,$\left\{\,p_{\,i}\,\right\}_{i \,=\, 1}^{\infty}$\; be a frame associated to \,$\left(\,a_{\,2},\, \cdots,\, a_{\,n}\,\right)$\, for \,$H$.\;Then the bounded linear operator \,$ T_{F} \,:\, l^{\,2}\,(\,\mathbb{N}\,) \,\to\, H_{F}$, defined by \,$T_{F}\,\{\,c_{i}\,\} \,=\, \sum\limits_{i \,=\, 1}^{\infty} \;c_{\,i}\,p_{\,i}$, is called  pre-frame operator and its adjoint operator described by
\[T_{F}^{\,\ast} \,:\, H_{F} \,\to\, l^{\,2}\,(\,\mathbb{N}\,),\;T_{F}^{\,\ast}\,p \,=\, \left \{\, \left <\,p,\, p_{i} \,|\, a_{\,2},\, \cdots,\, a_{\,n}\,\right >\,\right \}_{i \,=\, 1}^{\infty}\] is called the analysis operator.\;The operator \,$S_{F} \,:\, H_{F} \,\to\, H_{F}$\, given by 
\[S_{F}\,p \,=\, T_{F}\,T_{F}^{\,\ast}\,p \,=\, \sum\limits^{\infty}_{i \,=\, 1}\; \left <\,p \,,\, p_{\,i} \,|\, a_{\,2},\, \cdots,\, a_{\,n} \, \right >\,p_{\,i},\; \;\text{for all}\; \,p \,\in\, H_{F},\] is called the frame operator.
\end{definition}

\begin{theorem}\label{th2.21}
Let \,$\left\{\,p_{\,i}\,\right\}^{\,\infty}_{\,i \,=\, 1}$\; be a frame associated to \,$\left(\,a_{\,2},\, \cdots,\, a_{\,n}\,\right)$\; for \,$H$\, with bounds \,$A,\,B$.\;Then the corresponding frame operator \,$S_{F}$\, is bounded, invertible, self-adjoint and positive.
\end{theorem}

\begin{proof}
For each \,$p \,\in\, H_{F}$, we have 
\begin{align*}
&\left\|\,S_{F}\,p \,\right\|_{F}^{\,2} \,=\, \left\|\,S_{F}\,p,\, a_{\,2},\, \cdots,\,a_{\,n}\,\,\right\|^{\,2}\\
&=\, \sup\,\left\{\,\left|\,\left<\,S_{F}\,p,\, q \,|\, a_{\,2},\, \cdots,\, a_{\,n}\,\right>\,\right|^{\,2} \,:\, \left\|\,q,\, a_{\,2},\, \cdots,\, a_{\,n}\,\right\| \,=\, 1 \,\right\}\\
&=\, \sup\limits_{\left\|\,q,\, a_{\,2},\, \cdots,\, a_{\,n}\,\right\| \,=\, 1}\,\left|\,\left<\,\sum\limits^{\infty}_{i \,=\, 1}\,\left<\,p,\, p_{\,i} \,|\, a_{\,2},\, \cdots,\, a_{\,n} \,\right>\,p_{\,i},\, q \,|\, a_{\,2},\, \cdots,\, a_{\,n}\,\right>\,\right|^{\,2}\hspace{1cm}\\
&\leq\, \sup\limits_{\left\|\,q,\, a_{2},\, \cdots,\, a_{n}\,\right\| \,=\, 1}\,\sum\limits^{\infty}_{i \,=\, 1}\, \left|\, \left<\,p,\, p_{\,i} \,|\, a_{\,2},\, \cdots,\, a_{\,n} \,\right>\,\right|^{\,2}\,\sum\limits^{\infty}_{i \,=\, 1}\, \left|\,\left<\,q,\, p_{\,i} \,|\, a_{\,2},\, \cdots,\, a_{\,n}\,\right>\,\right|^{\,2}\\
&[\;\text{using Cauchy-Schwarz iequality}\;]\\
&\leq\, B^{\,2}\, \left \|\,p,\, a_{\,2},\, \cdots,\, a_{\,n}\,\right\|^{\,2} \,=\, B^{\,2}\, \|\,p\,\|^{\,2}_{F}\hspace{5.5cm}\\
&[\,\text{since}\, \left\{\,p_{\,i}\,\right\}^{\,\infty}_{\,i \,=\, 1}\; \;\text{is a frame associated to \,$\left(\,a_{\,2},\, \cdots,\, a_{\,n}\,\right)$}\;\;]
\end{align*}
This shows that \,$S_{F}$\; is bounded.\;Since \,$S_{F} \,=\, T_{F}\,T^{\,\ast}_{F}$, it is easy to verify that \,$S_{F}$\, is self-adjoint.\,The inequality (\ref{eqp1}), can be written as 
\[A\left<\,p,\, p \,|\, a_{\,2},\, \cdots,\, a_{\,n}\,\right> \,\leq \left<\,S_{F}\,p,\, p \,|\, a_{\,2},\, \cdots,\, a_{\,n}\right> \leq B\left<\,p,\, p \,|\, a_{\,2},\, \cdots,\, a_{\,n}\,\right>\]
and this gives \,$A\,I_{F} \,\leq\, S_{F} \,\leq\, B\,I_{F}$.\,Thus, \,$S_{F}$\, is positive and consequently it is invertible.    
\end{proof}

\begin{remark}\label{rmr1}
In Theorem \ref{th2.21}, it is proved that \,$A\,I_{F} \,\leq\, S_{F} \,\leq\, B\,I_{F}$. Since \,$S^{\,-1}_{F}$\; commutes with both \,$S_{F}$\; and \,$I_{F}$, multiplying in the inequality, \,$ A\,I_{F} \,\leq\, S_{F} \,\leq\, B\,I_{F}$\, by \,$S^{\,-1}_{F}$, we get \,$B^{\,-1}\,I_{F} \,\leq\, S^{\,-1}_{F} \,\leq\, A^{\,-1}\,I_{F}$.
\end{remark}

For more details on frame in \,$n$-Hilbert space one can go through the paper \cite{Prasenjit}. \\

Let \,$H$\, and \,$K$\, be two \,$n$-Hilbert spaces associated with the \,$n$-inner products \,$\left<\,\cdot,\, \cdot \,|\, \cdot,\, \cdots,\, \cdot\,\right>_{1}$\, and \,$\left<\,\cdot,\, \cdot \,|\, \cdot,\, \cdots,\, \cdot\,\right>_{2}$, respectively.\;The tensor product of \,$H$\, and \,$K$\, is denoted by \,$H \,\otimes\, K$\, and it is defined to be an \,$n$-inner product space associated with the \,$n$-inner product given by 
\[\left<\,f_{\,1} \,\otimes\, g_{\,1},\, f_{\,2} \,\otimes\, g_{\,2} \,|\, f_{\,3} \,\otimes\, g_{\,3},\, \,\cdots,\, f_{\,n} \,\otimes\, g_{\,n}\,\right>\]
\begin{equation}\label{eqn1}
 \,=\, \left<\,f_{\,1},\, f_{\,2} \,|\, f_{\,3},\, \,\cdots,\, f_{\,n}\,\right>_{1}\,\left<\,g_{\,1},\, g_{\,2} \,|\, g_{\,3},\, \,\cdots,\, g_{\,n}\,\right>_{2},
\end{equation}
for all \,$f_{\,1},\, f_{\,2},\, f_{\,3},\, \,\cdots,\, f_{\,n} \,\in\, H$\, and \,$g_{\,1},\, g_{\,2},\, g_{\,3},\, \,\cdots,\, g_{\,n} \,\in\, K$.\\
The \,$n$-norm on \,$H \,\otimes\, K$\, is defined by 
\[\left\|\,f_{\,1} \,\otimes\, g_{\,1},\, f_{\,2} \,\otimes\, g_{\,2},\, \,\cdots,\,\, f_{\,n} \,\otimes\, g_{\,n}\,\right\|\]
\begin{equation}\label{eqn1.1}
\hspace{.6cm} =\,\left\|\,f_{\,1},\, f_{\,2},\, \cdots,\, f_{\,n}\,\right\|_{1}\;\left\|\,g_{\,1},\, g_{\,2},\, \cdots,\, g_{\,n}\,\right\|_{2},
\end{equation}
for all \,$f_{\,1},\, f_{\,2},\, \,\cdots,\, f_{\,n} \,\in\, H\, \;\text{and}\; \,g_{\,1},\, g_{\,2},\, \,\cdots,\, g_{\,n} \,\in\, K$, where the \,$n$-norms \,$\left\|\,\cdot,\, \cdots,\, \cdot \,\right\|_{1}$\, and \,$\left\|\,\cdot,\, \cdots,\, \cdot \,\right\|_{2}$\, are generated by \,$\left<\,\cdot,\, \cdot \,|\, \cdot,\, \cdots,\, \cdot\,\right>_{1}$\, and \,$\left<\,\cdot,\, \cdot \,|\, \cdot,\, \cdots,\, \cdot\,\right>_{2}$, respectively.\;The space \,$H \,\otimes\, K$\, is complete with respect to the above \,$n$-inner product.\;Therefore the space \,$H \,\otimes\, K$\, is an \,$n$-Hilbert space.

Consider \,$G \,=\, \left\{\,b_{\,2},\, b_{\,3},\, \cdots,\, b_{\,n}\,\right\}$, where \,$b_{\,2},\, b_{\,3},\, \cdots,\, b_{\,n}$\, are fixed elements in \,$K$\, and \,$L_{G}$\, denote the linear subspace of \,$K$\, spanned by \,$G$.\,Now, we can define the Hilbert space \,$K_{G}$\, with respect to the inner product is given by
\[\left<\,p \,+\, L_{G},\, q \,+\, L_{G}\,\right>_{G} \,=\, \left<\,p,\, q\,\right>_{G} \,=\, \left<\,p,\, q \,|\, b_{\,2},\,  \cdots,\, b_{\,n} \,\right>_{2}; \;\forall \;\; p,\, q \,\in\, K.\]

\begin{remark}
According to the definition \ref{def0.001}, \,$H_{F} \,\otimes\, K_{G}$\, is the Hilbert space with respect to the inner product:
\[\left<\,p \,\otimes\, q,\, p^{\,\prime} \,\otimes\, q^{\,\prime}\,\right> \,=\, \left<\,p,\, p^{\,\prime}\,\right>_{F}\;\left<\,q,\, q^{\,\prime}\,\right>_{G},\]
for all \,$p,\, p^{\,\prime} \,\in\, H_{F}\; \;\text{and}\; \;q,\, q^{\,\prime} \,\in\, K_{G}$.    
\end{remark}

\begin{remark}
From the definition of ordinary frames for separable Hilbert spaces, the sequence of vectors \,$\left\{\,p_{\,i} \,\otimes\, q_{\,j}\,\right\}^{\,\infty}_{i,\,j \,=\, 1}$\,  in \,$H \,\otimes\, K$\, can be consider as a frame associated to \,$\left(\,a_{\,2} \,\otimes\, b_{\,2},\, \,\cdots,\, a_{\,n} \,\otimes\, b_{\,n}\,\right)$\, for \,$H \,\otimes\, K$\, if there exist constants \,$0 \,<\, A \,\leq\, B \,<\, \infty$\, such that
\begin{align}\label{eqp1.1}
&A \left\|\,p \,\otimes\, q,\, a_{\,2} \,\otimes\, b_{\,2},\,\cdots,\, a_{\,n} \,\otimes\, b_{\,n}\,\right\|^{\,2}\\
& \leq \sum\limits_{i,\, j \,=\, 1}^{\,\infty}\left|\,\left<\,p \,\otimes\, q,\, p_{\,i} \,\otimes\, q_{\,j} \,|\, a_{\,2} \,\otimes\, b_{\,2},\, \cdots,\, a_{\,n} \,\otimes\, b_{\,n}\,\right>\,\right|^{\,2}\nonumber\\
&\leq\, B\,\left\|\,p \,\otimes\, q,\, a_{\,2} \,\otimes\, b_{\,2},\, \cdots,\, a_{\,n} \,\otimes\, b_{\,n}\,\right\|^{\,2}; \;\forall\, \;p \,\otimes\, q \,\in\, H \,\otimes\, K,\nonumber
\end{align} 
where \,$\{\,p_{\,i}\,\}_{i \,=\,1}^{\infty}$\, and \,$\{\,q_{\,j}\,\}_{j \,=\,1}^{\infty}$\, be the sequences of vectors in \,$H$\, and \,$K$, respectively and \,$a_{\,2} \,\otimes\, b_{\,2},\, a_{\,3} \,\otimes\, b_{\,3},\,\cdots,\, a_{\,n} \,\otimes\, b_{\,n}$\, be fixed elements in \,$H \,\otimes\, K$.\;The constants \,$A,\,B$\, are called the frame bounds.\;If \,$A \,=\, B$\, then it is called a tight frame associated to \,$(\,a_{\,2} \,\otimes\, b_{\,2},\, \,\cdots,\, a_{\,n} \,\otimes\, b_{\,n}\,)$.\;If the sequence \,$\left\{\,p_{\,i} \,\otimes\, q_{\,j}\,\right\}^{\,\infty}_{i,\,j \,=\, 1}$\, satisfies the right hand inequality of (\ref{eqp1.1}), is called a Bessel sequence associated to \,$\left(\,a_{\,2} \,\otimes\, b_{\,2},\, \,\cdots,\, a_{\,n} \,\otimes\, b_{\,n}\,\right)$\, in \,$H \,\otimes\, K$.    
\end{remark}

\begin{theorem}\label{th2.1}
Let \,$\{\,p_{\,i}\,\}_{i \,=\,1}^{\infty}$\, and \,$\{\,q_{\,j}\,\}_{j \,=\,1}^{\infty}$\, be the sequences of vectors in \,$n$-Hilbert spaces \,$H$\, and \,$K$.\;The sequence \,$\left\{\,p_{\,i} \,\otimes\, q_{\,j}\,\right\}^{\,\infty}_{i,\,j \,=\, 1} \,\subseteq\, H \,\otimes\, K$\, is a frame associated to \,$\left(\,a_{\,2} \,\otimes\, b_{\,2},\, \,\cdots,\, a_{\,n} \,\otimes\, b_{\,n}\,\right)$\, for \,$H \,\otimes\, K$\, if and only if \,$\{\,p_{\,i}\,\}_{i \,=\,1}^{\infty}$\, is a frame associated to \,$\left(\,a_{\,2},\, \cdots,\, a_{\,n}\,\right)$\, for \,$H$\, and \,$\{\,q_{\,j}\,\}_{j \,=\,1}^{\infty}$\, is a frame associated to \,$\left(\,b_{\,2},\, \cdots,\, b_{\,n}\,\right)$\, for \,$K$.   
\end{theorem}

\begin{proof}
Suppose that the sequence \,$\left\{\,p_{\,i} \,\otimes\, q_{\,j}\,\right\}^{\,\infty}_{i,\,j \,=\, 1}$\, is a frame associated to \,$\left(\,a_{\,2} \,\otimes\, b_{\,2},\, \,\cdots,\, a_{\,n} \,\otimes\, b_{\,n}\,\right)$\, for \,$H \,\otimes\, K$.\;Then, for each \,$p \,\otimes\, q \,\in\, H \,\otimes\, K \,-\, \{\,\theta \,\otimes\, \theta\,\}$, there exist constants \,$A,\,B \,>\, 0$\, such that
\begin{align*}
&A \left\|\,p \,\otimes\, q,\, a_{\,2} \,\otimes\, b_{\,2},\,\cdots,\, a_{\,n} \,\otimes\, b_{\,n}\,\right\|^{\,2}\hspace{2.7cm}\\
&\leq \sum\limits_{i,\, j \,=\, 1}^{\,\infty}\left|\,\left<\,p \,\otimes\, q,\, p_{\,i} \,\otimes\, q_{\,j} \,|\, a_{\,2} \,\otimes\, b_{\,2},\, \cdots,\, a_{\,n} \,\otimes\, b_{\,n}\,\right>\,\right|^{\,2}\\
&\leq\, B\,\left\|\,p \,\otimes\, q,\, a_{\,2} \,\otimes\, b_{\,2},\, \cdots,\, a_{\,n} \,\otimes\, b_{\,n}\,\right\|^{\,2}.
\end{align*}
Using the \,$n$-norm and \,$n$-inner product on \,$H \,\otimes\, K$, we get
\begin{align*}
&A\left\|\,p,\, a_{2},\, \cdots,\, a_{n}\,\right\|_{1}^{\,2}\left\|\,q,\, b_{2},\, \cdots,\, b_{n}\,\right\|_{2}^{\,2}\leq \left(\sum\limits^{\infty}_{i \,=\, 1}\left|\,\left<\,p,\, p_{\,i}\,|\,a_{\,2},\, \cdots,\, a_{\,n}\,\right>_{1}\,\right|^{\,2}\right)\times\\
&\left(\sum\limits_{\,j \,=\, 1}^{\,\infty}\left|\,\left<\,q,\, q_{\,j} \,|\, b_{2},\, \cdots,\, b_{n}\,\right>_{2}\,\right|^{\,2}\right)  \leq B\left\|\,p,\, a_{2},\, \cdots,\, a_{n}\,\right\|_{1}^{\,2}\,\left\|\,q,\, b_{2},\, \cdots,\, b_{n}\,\right\|_{2}^{\,2}.
\end{align*}
Since \,$p \,\otimes\, q \,\in\, H \,\otimes\, K$\, is  non-zero element i.\,e., \,$p \,\in\, H$\, and \,$q \,\in\, K$\, are non-zero elements.\,Here, we may assume that every \,$p_{\,i}$\, and \,$a_{\,2},\, \cdots$, \,$ a_{\,n}$\, are linearly independent and every \,$q_{\,j}$\, and \,$b_{\,2},\, \cdots$, \,$ b_{\,n}$\, are linearly independent.\,Hence 
\[\sum\limits_{\,j \,=\, 1}^{\,\infty}\,\left|\,\left<\,q,\, q_{\,j} \,|\, b_{\,2},\, \cdots,\, b_{\,n}\,\right>_{2}\,\right|^{\,2},\,\sum\limits_{\,i \,=\, 1}^{\,\infty}\,\left|\,\left<\,p ,\, p_{\,i} \,|\, a_{\,2},\, \cdots,\, a_{\,n}\,\right>_{1}\,\right|^{\,2}\] are non-zero.\,Therefore, by the above inequality, we get   
\begin{align*}
&\dfrac{A\,\left\|\,q,\, b_{\,2},\, \cdots,\, b_{\,n}\,\right\|_{2}^{\,2}\,\left\|\,p,\, a_{\,2},\, \cdots,\, a_{\,n}\,\right\|_{1}^{\,2}}{\sum\limits_{\,j \,=\, 1}^{\,\infty}\,\left|\,\left<\,q,\, q_{\,j} \,|\, b_{\,2},\, \cdots,\, b_{\,n}\,\right>_{2}\,\right|^{\,2}}\,\leq\, \sum\limits_{\,i \,=\, 1}^{\,\infty}\,\left|\,\left<\,p,\, p_{\,i} \,|\, a_{\,2},\, \cdots,\, a_{\,n}\,\right>_{1}\,\right|^{\,2}\\
&\leq\, \dfrac{B\,\left\|\,q,\, b_{\,2},\, \cdots,\, b_{\,n}\,\right\|_{2}^{\,2}}{\sum\limits_{\,j \,=\, 1}^{\,\infty}\,\left|\,\left<\,q,\, q_{\,j} \,|\, b_{\,2},\, \cdots,\, b_{\,n}\,\right>_{2}\,\right|^{\,2}}\,\left\|\,p,\, a_{\,2},\, \cdots,\, a_{\,n}\,\right\|_{1}^{\,2}.
\end{align*}
This implies that
\[A_{1}\left\|\,p,\, a_{2},\, \cdots,\, a_{n}\,\right\|_{1}^{\,2} \leq \sum\limits_{\,i \,=\, 1}^{\,\infty}\left|\,\left<\,p,\, p_{\,i}\,|\,a_{2},\, \cdots,\, a_{n}\,\right>_{1}\,\right|^{\,2} \leq\, B_{1}\left\|\,p,\, a_{2},\, \cdots,\, a_{n}\,\right\|_{1}^{\,2},\]
for all \,$p \,\in\, H$, where \,$A_{1} \,=\, \inf\limits_{q \,\in\, K}\left\{\,\dfrac{A\,\left\|\,q,\, b_{2},\, \cdots,\, b_{n}\,\right\|_{2}^{\,2}}{\sum\limits_{\,j \,=\, 1}^{\,\infty}\,\left|\,\left<\,q,\, q_{\,j} \,|\, b_{2},\, \cdots,\, b_{n}\,\right>_{2}\,\right|^{\,2}}\,\right\}$\\ and \,$B_{1} \,=\, \sup\limits_{q \,\in\, K}\left\{\,\dfrac{B\,\left\|\,q,\, b_{\,2},\, \cdots,\, b_{\,n}\,\right\|_{2}^{\,2}}{\sum\limits_{\,j \,=\, 1}^{\,\infty}\,\left|\,\left<\,q,\, q_{\,j} \,|\, b_{\,2},\, \cdots,\, b_{\,n}\,\right>_{2}\,\right|^{\,2}}\,\right\}$.\\This shows that \,$\{\,p_{\,i}\,\}_{i \,=\,1}^{\infty}$\, is a frame associated to \,$\left(\,a_{\,2},\, \cdots,\, a_{\,n}\,\right)$\, for \,$H$.\;Similarly, it can be shown that \,$\{\,q_{\,j}\,\}_{j \,=\,1}^{\infty}$\, is a frame associated to \,$\left(\,b_{\,2},\, \cdots,\, b_{\,n}\,\right)$\, for \,$K$.\\

Conversely, suppose that \,$\{\,p_{\,i}\,\}_{i \,=\,1}^{\infty}$\, is a frame associated to \,$\left(\,a_{2},\, \cdots,\, a_{n}\,\right)$\, for \,$H$\, with bounds \,$A,\,B$\, and \,$\{\,q_{\,j}\,\}_{j \,=\,1}^{\infty}$\, is a frame associated to \,$\left(\,b_{2},\, \cdots,\, b_{n}\,\right)$\, for \,$K$\, with bounds \,$C,\,D$.\;Then, for all \,$p \,\in\, H$\, and \,$q \,\in\, K$, we have
\begin{align*}
&A\left\|\,p,\, a_{2},\, \cdots,\, a_{n}\,\right\|_{1}^{\,2} \leq \sum\limits_{\,i \,=\, 1}^{\,\infty}\left|\,\left<\,p,\, p_{\,i}\,|\,a_{2},\, \cdots,\, a_{n}\,\right>_{1}\,\right|^{\,2} \leq B\left\|\,p,\, a_{2},\, \cdots,\, a_{n}\,\right\|_{1}^{\,2},\\
&C\left\|\,q,\, b_{2},\, \cdots,\, b_{n}\,\right\|_{2}^{\,2} \leq \sum\limits_{\,j \,=\, 1}^{\,\infty}\left|\,\left<\,q,\, q_{\,j}\,|\,b_{2},\, \cdots,\, b_{n}\,\right>_{2}\,\right|^{\,2} \leq D\left\|\,q,\, b_{2},\, \cdots,\, b_{n}\,\right\|_{2}^{\,2}.
\end{align*}
Multiplying the above two inequalities and using (\ref{eqn1}) and (\ref{eqn1.1}), we get
\begin{align*}
&A\,C \left\|\,p \,\otimes\, q,\, a_{\,2} \,\otimes\, b_{\,2},\,\cdots,\, a_{\,n} \,\otimes\, b_{\,n}\,\right\|^{\,2}\hspace{2.6cm}\\
& \leq \sum\limits_{i,\, j \,=\, 1}^{\,\infty}\left|\,\left<\,p \,\otimes\, q,\, p_{\,i} \,\otimes\, q_{\,j} \,|\, a_{\,2} \,\otimes\, b_{\,2},\, \cdots,\, a_{\,n} \,\otimes\, b_{\,n}\,\right>\,\right|^{\,2}\\
&\leq\, B\,D\left\|\,p \,\otimes\, q,\, a_{\,2} \,\otimes\, b_{\,2},\, \cdots,\, a_{\,n} \,\otimes\, b_{\,n}\,\right\|^{\,2}.
\end{align*}
Hence, \,$\left\{\,p_{\,i} \otimes q_{\,j}\,\right\}^{\,\infty}_{i,\,j \,=\, 1}$\, is a frame associated to \,$\left(\,a_{\,2} \otimes b_{\,2},\, \,\cdots,\, a_{\,n} \otimes b_{\,n}\,\right)$\, for \,$H \,\otimes\, K$. 
\end{proof}

\begin{theorem}
The sequence \,$\left\{\,p_{\,i} \,\otimes\, q_{\,j}\,\right\}^{\,\infty}_{i,\,j \,=\, 1} \,\subseteq\, H \,\otimes\, K$\, is a frame associated to \,$\left(\,a_{\,2} \,\otimes\, b_{\,2},\, \,\cdots,\, a_{\,n} \,\otimes\, b_{\,n}\,\right)$\, for \,$H \,\otimes\, K$\, if and only if \,$\left\{\,p_{\,i} \,\otimes\, q_{\,j}\,\right\}^{\,\infty}_{i,\,j \,=\, 1}$\, is a frame  for \,$H_{F} \,\otimes\, K_{G}$. 
\end{theorem}

\begin{proof}
Suppose that the sequence \,$\left\{\,p_{\,i} \,\otimes\, q_{\,j}\,\right\}^{\,\infty}_{i,\,j \,=\, 1}$\, is a frame associated to \,$\left(\,a_{\,2} \,\otimes\, b_{\,2},\, \,\cdots,\, a_{\,n} \,\otimes\, b_{\,n}\,\right)$\, for \,$H \,\otimes\, K$.\;By Theorem \ref{th2.1}, \,$\{\,p_{\,i}\,\}_{i \,=\,1}^{\infty}$\, is a frame associated to \,$\left(\,a_{\,2},\, \cdots,\, a_{\,n}\,\right)$\, for \,$H$\, and \,$\{\,q_{\,j}\,\}_{j \,=\,1}^{\infty}$\, is a frame associated to \,$\left(\,b_{\,2},\, \cdots,\, b_{\,n}\,\right)$\, for \,$K$.\;Now, applying the Theorem \ref{th2}, \,$\{\,p_{\,i}\,\}_{i \,=\,1}^{\infty}$\, and \,$\{\,q_{\,j}\,\}_{j \,=\,1}^{\infty}$\, are frames for \,$H_{F}$\, and \,$K_{G}$, respectively.\;Hence, \,$\left\{\,p_{\,i} \,\otimes\, q_{\,j}\,\right\}^{\,\infty}_{i,\,j \,=\, 1}$\, is a frame  for \,$H_{F} \,\otimes\, K_{G}$.\\The proof of the converse part is obvious.      
\end{proof}

\begin{remark}
Let \,$\left\{\,p_{\,i} \,\otimes\, q_{\,j}\,\right\}^{\,\infty}_{i,\,j \,=\, 1}$\, be a frame associated to \,$(\,a_{\,2} \,\otimes\, b_{\,2}$, \,$\cdots,\, a_{\,n} \,\otimes\, b_{\,n}\,)$\, for \,$H \,\otimes\, K$.\;According to the definition \ref{def0.0001}, the frame operator \,$S_{F \,\otimes\, G} \,:\, H_{F} \,\otimes\, K_{G} \,\to\, H_{F} \,\otimes\, K_{G}$\, is described by 
\[S_{F \otimes G}\,(\,p \,\otimes\, q\,) = \sum\limits_{i,\, j \,=\, 1}^{\,\infty}\left<\,p \otimes q,\, p_{\,i} \otimes q_{\,j} \,|\, a_{\,2} \otimes b_{\,2}, \,\cdots,\, a_{\,n} \otimes b_{\,n}\,\right>\,\left(\,p_{\,i} \otimes q_{\,j}\,\right)\]
for all \,$p \,\otimes\, q \,\in\, H_{F} \,\otimes\, K_{G}$. 
\end{remark}

\begin{theorem}
If \,$S_{F},\, \,S_{G}$\, and \,$S_{F \,\otimes\, G}$\, are the corresponding frame operator for \,$\{\,p_{\,i}\,\}_{i \,=\,1}^{\infty},\, \,\{\,q_{\,j}\,\}_{j \,=\,1}^{\infty}$\, and \,$\left\{\,p_{\,i} \,\otimes\, q_{\,j}\,\right\}^{\,\infty}_{i,\,j \,=\, 1}$, respectively, then \,$S_{F \,\otimes\, G} \,=\, S_{F} \,\otimes\, S_{G}$\, and \,$S^{\,-\, 1}_{F \,\otimes\, G} \,=\, S^{\,-\, 1}_{F} \,\otimes\, S^{\,-\, 1}_{G}$.  
\end{theorem}

\begin{proof}
Since \,$S_{F \,\otimes\, G}$\, is the frame operator for \,$\left\{\,p_{\,i} \,\otimes\, q_{\,j}\,\right\}^{\,\infty}_{i,\,j \,=\, 1}$, we have
\begin{align*}
&S_{F \otimes G}\,(\,p \,\otimes\, q\,) = \sum\limits_{i,\, j \,=\, 1}^{\,\infty}\left<\,p \otimes q,\, p_{\,i} \otimes q_{\,j} \,|\, a_{\,2} \otimes b_{\,2}, \,\cdots,\, a_{\,n} \otimes b_{\,n}\,\right>\,\left(\,p_{\,i} \otimes q_{\,j}\,\right)\\
&\,=\, \sum\limits_{i,\, j \,=\, 1}^{\,\infty}\,\left<\,p,\, p_{\,i} \,|\, a_{\,2},\, \cdots,\, a_{\,n}\,\right>_{1}\,\left<\,q,\, q_{\,j} \,|\, b_{\,2},\, \cdots,\, b_{\,n}\,\right>_{2}\,\left(\,p_{\,i} \,\otimes\, q_{\,j}\,\right)\\
&= \left(\sum\limits_{\,i \,=\, 1}^{\,\infty}\left<\,p,\, p_{\,i} \,|\, a_{\,2},\, \cdots,\, a_{\,n}\,\right>_{1}\,p_{\,i}\right) \,\otimes\, \left(\sum\limits_{\,j \,=\, 1}^{\,\infty}\,\left<\,q,\, q_{\,j} \,|\, b_{\,2},\, \cdots,\, b_{\,n}\,\right>_{2}\,q_{\,j}\right)\\
&=\, S_{F}\,(\,p\,) \,\otimes\, S_{G}\,(\,q\,) \,=\, \left(\,S_{F} \,\otimes\, S_{G}\,\right)\,(\,p \,\otimes\, q\,)\; \;\forall\; p \,\otimes\, q \,\in\, H_{F} \,\otimes\, K_{G}.
\end{align*}
Thus \,$S_{F \,\otimes\, G} \,=\, S_{F} \,\otimes\, S_{G}$.\;Since \,$S_{F}\, \,\text{and}\,S_{G}$\, are invertible, by Theorem \ref{th1.1} \,$(iv)$, \,$S^{\,-\, 1}_{F \,\otimes\, G} \,=\, \left(\,S_{F} \,\otimes\, S_{G}\,\right)^{\,-\, 1} \,=\, S^{\,-\, 1}_{F} \,\otimes\, S^{\,-\, 1}_{G}$.    
\end{proof}

\begin{theorem}\label{th2.2}
Let \,$\{\,p_{\,i}\,\}_{i \,=\,1}^{\infty}$\, be a frame associated to \,$\left(\,a_{\,2},\, \cdots,\, a_{\,n}\,\right)$\, for \,$H$\, with bounds \,$A,\,B$\, and \,$\{\,q_{\,j}\,\}_{j \,=\,1}^{\infty}$\, be a frame associated to \,$(\,b_{\,2},\, \cdots,$\, \,$b_{\,n}\,)$\, for \,$K$\, with bounds \,$C,\,D$\, having their corresponding frame operators \,$S_{F}$\, and \,$S_{G}$, respectively.\;Then \,$\Lambda \,=\, \left\{\,S^{\,-\, 1}_{F \,\otimes\, G}\,\left(\,p_{\,i} \,\otimes\, q_{\,j}\,\right)\,\right\}^{\infty}_{i,\,j \,=\, 1}$ is a frame associated to \,$\left(a_{2} \otimes b_{2},\, \,\cdots,\, a_{n} \otimes b_{n}\right)$\, for \,$H \,\otimes\, K$\, with the corresponding frame operator \,$S^{\,-\, 1}_{F \,\otimes\, G}$.    
\end{theorem}

\begin{proof}
For each \,$p \,\otimes\, q \,\in\, H_{F} \,\otimes\, K_{G}$, we have
\begin{align}
&\sum\limits_{i,\, j \,=\, 1}^{\,\infty}\,\left|\,\left<\,p \,\otimes\, q,\, S^{\,-\, 1}_{F \,\otimes\, G}\,\left(\,p_{\,i} \,\otimes\, q_{\,j}\,\right) \,|\, a_{\,2} \,\otimes\, b_{\,2}, \,\cdots,\, a_{\,n} \,\otimes\, b_{\,n}\,\right>\,\right|^{\,2}\nonumber\\
&= \sum\limits_{i,\, j \,=\, 1}^{\,\infty}\left|\,\left<\,p \,\otimes\, q,\, \left(\,S^{\,-\, 1}_{F} \otimes S^{\,-\, 1}_{G}\,\right)\,\left(\,p_{\,i} \otimes q_{\,j}\,\right) \,|\, a_{\,2} \,\otimes\, b_{\,2}, \,\cdots,\, a_{\,n} \,\otimes\, b_{\,n}\,\right>\,\right|^{\,2}\nonumber\\
&=\, \sum\limits_{i,\, j \,=\, 1}^{\,\infty}\,\left|\,\left<\,p \,\otimes\, q,\,S^{\,-\, 1}_{F}\,\left(\,p_{\,i}\,\right) \,\otimes\, S^{\,-\, 1}_{G}\,\left(\,q_{\,j}\,\right) \,|\, a_{\,2} \,\otimes\, b_{\,2}, \,\cdots,\, a_{\,n} \,\otimes\, b_{\,n}\,\right>\,\right|^{\,2}\nonumber\\
&= \sum\limits_{\,i \,=\, 1}^{\,\infty}\left|\,\left<\,p,\, S^{\,-\, 1}_{F}\,p_{\,i}\,|\,a_{\,2},\, \cdots,\, a_{\,n}\,\right>_{1}\,\right|^{\,2}\sum\limits_{\,j \,=\, 1}^{\,\infty}\left|\,\left<\,q,\, S^{\,-\, 1}_{G}\,q_{\,j}\,|\, b_{\,2},\, \cdots,\, b_{\,n}\,\right>_{2}\,\right|^{\,2}\nonumber\\
&=\, \sum\limits_{\,i \,=\, 1}^{\,\infty}\,\left|\,\left<\,S^{\,-\, 1}_{F}\,p,\, p_{\,i}\,|\,a_{\,2},\, \cdots,\, a_{\,n}\,\right>_{1}\,\right|^{\,2}\sum\limits_{\,j \,=\, 1}^{\,\infty}\left|\,\left<\,S^{\,-\, 1}_{G}\,q,\, q_{\,j}\,|\,b_{\,2},\, \cdots,\, b_{\,n}\,\right>_{2}\,\right|^{\,2}\label{eq1.2}\\
&\leq\, B\,\left\|\,S^{\,-\, 1}_{F}\,(\,p\,),\, a_{\,2},\, \cdots,\, a_{\,n}\,\right\|_{1}^{\,2}\,D\,\left\|\,S^{\,-\, 1}_{G}\,(\,q\,),\, b_{\,2},\, \cdots,\, b_{\,n}\,\right\|_{2}^{\,2}\nonumber\\
&[\;\text{since \,$\{\,p_{\,i}\,\}_{i \,=\,1}^{\infty}$\, is a frame associated to \,$\left(\,a_{\,2},\, \cdots,\, a_{\,n}\,\right)$, and}\nonumber\\
&\text{\,$\{\,q_{\,j}\,\}_{j \,=\,1}^{\infty}$\, is a frame associated to \,$\left(\,b_{\,2},\, \cdots,\, b_{\,n}\,\right)$\,}\;]\nonumber\\
&\leq\, B\,D\,\left\|\,S^{\,-\, 1}_{F}\,\right\|^{\,2}\,\left\|\,S^{\,-\, 1}_{G}\,\right\|^{\,2}\,\left\|\,p,\, a_{\,2},\, \cdots,\, a_{\,n}\,\right\|_{1}^{\,2}\,\left\|\,q,\, b_{\,2},\, \cdots,\, b_{\,n}\,\right\|_{2}^{\,2}\nonumber\\
&\left[\;\text{since $S^{\,-\, 1}_{F},\,S^{\,-\, 1}_{G}$ are bounded operators on $H_{F}$ and $K_{G}$, respectively}\;\right]\nonumber\\
&\leq\, \dfrac{B\,D}{A^{\,2}\,C^{\,2}}\,\left\|\,p \,\otimes\, q,\, a_{\,2} \,\otimes\, b_{\,2} \,\cdots,\, a_{\,n} \,\otimes\, b_{\,n}\,\right\|^{\,2}\;[\;\text{using (\ref{eqn1.1})}\;]\nonumber\\
&[\;\text{since $B^{\,-1}\,I_{F} \leq S^{\,-1}_{F} \leq A^{\,-1}\,I_{F}$\, and \,$D^{\,-1}\,I_{G} \leq S^{\,-1}_{G} \leq C^{\,-1}\,I_{G}$},\nonumber\\
&\text{where \,$I_{G}$\, denote the identity operator on \,$K_{G}$}\;].\nonumber
\end{align}
Now, for \,$p \,\in\, H_{F}$, we have
\[\left\|\,p,\, a_{\,2},\, \cdots,\, a_{\,n}\,\right\|_{1} \,\leq\, \left\|\,S_{F}\,\right\|\,\left\|\,S^{\,-\, 1}_{F}\,(\,p\,),\, a_{\,2},\, \cdots,\, a_{\,n}\,\right\|_{1}\]
\begin{equation}\label{eq1.3}
\hspace{2.65cm}\,\leq\, B\,\left\|\,S^{\,-\, 1}_{F}\,(\,p\,),\, a_{\,2} \,,\, \cdots \,,\, a_{\,n}\,\right\|_{1}
\end{equation}
and similarly for \,$q \,\in\, K_{G}$, we get 
\begin{equation}\label{eq1.4}
\left\|\,q,\, b_{\,2} \,,\, \cdots \,,\, b_{\,n}\,\right\|_{2} \,\leq\, D\,\left\|\,S^{\,-\, 1}_{G}\,(\,q\,),\, b_{\,2} \,,\, \cdots \,,\, b_{\,n}\,\right\|_{2}.
\end{equation} 
On the other hand, from (\ref{eq1.2}),
\[\sum\limits_{i,\, j \,=\, 1}^{\,\infty}\,\left|\,\left<\,p \,\otimes\, q,\, S^{\,-\, 1}_{F \,\otimes\, G}\,\left(\,p_{\,i} \,\otimes\, q_{\,j}\,\right) \,|\, a_{\,2} \,\otimes\, b_{\,2}, \,\cdots,\, a_{\,n} \,\otimes\, b_{\,n}\,\right>\,\right|^{\,2}\hspace{3cm}\]
\[= \sum\limits_{\,i \,=\, 1}^{\,\infty}\left|\,\left<\,S^{\,-\, 1}_{F}p,\, p_{\,i}\,|\,a_{\,2},\, \cdots,\, a_{\,n}\,\right>_{1}\,\right|^{\,2}\sum\limits_{\,j \,=\, 1}^{\,\infty}\left|\,\left<\,S^{\,-\, 1}_{G}q,\, q_{\,j}\,|\,b_{\,2},\, \cdots,\, b_{\,n}\,\right>_{2}\,\right|^{\,2}\]
\[\geq\,  A\,\left\|\,S^{\,-\, 1}_{F}\,(\,p\,),\, a_{\,2} \,,\, \cdots \,,\, a_{\,n}\,\right\|_{1}^{\,2}\,C\,\left\|\,S^{\,-\, 1}_{G}\,(\,q\,),\, b_{\,2} \,,\, \cdots \,,\, b_{\,n}\,\right\|_{2}^{\,2}\hspace{3cm}\]
\[\geq\, \dfrac{A\,C}{B^{\,2}\,D^{\,2}}\,\left\|\,p,\, a_{\,2},\, \cdots,\, a_{\,n}\,\right\|_{1}^{\,2}\,\left\|\,q,\, b_{\,2},\, \cdots,\, b_{\,n}\,\right\|_{2}^{\,2}\;[\;\text{by (\ref{eq1.3}) and (\ref{eq1.4})}\;]\hspace{2.2cm}\]
\[ \,=\, \dfrac{A\,C}{B^{\,2}\,D^{\,2}}\,\left\|\,p \,\otimes\, q,\, a_{\,2} \,\otimes\, b_{\,2} \,\cdots,\, a_{\,n} \,\otimes\, b_{\,n}\,\right\|^{\,2}.\hspace{4.3cm}\]
Hence, \,$\Lambda$\, is a frame associated to \,$\left(\,a_{\,2} \,\otimes\, b_{\,2},\, \,\cdots,\, a_{\,n} \,\otimes\, b_{\,n}\,\right)$\, for \,$H \,\otimes\, K$\, with bounds \,$\dfrac{A\,C}{B^{\,2}\,D^{\,2}}$\, and \,$\dfrac{B\,D}{A^{\,2}\,C^{\,2}}$.\\

Furthermore, for each \,$p \,\otimes\, q \,\in\, H_{F} \,\otimes\, K_{G}$, we have
\begin{align*}
&\sum\limits_{i,\, j \,=\, 1}^{\,\infty}\left<\,p \otimes q,\, S^{\,-\, 1}_{F \otimes G}\left(p_{\,i} \otimes q_{\,j}\right)\,|\,a_{\,2} \otimes b_{\,2}, \,\cdots,\, a_{\,n} \otimes b_{\,n}\,\right>S^{\,-\, 1}_{F \otimes G}\,\left(p_{\,i} \otimes q_{\,j}\right)\\
&= \sum\limits_{i,\, j \,=\, 1}^{\,\infty}\left<\,p \otimes q,\, S^{\,-\, 1}_{F}p_{i} \otimes S^{\,-\, 1}_{G}q_{j}\,|\,a_{2} \otimes b_{2}, \,\cdots,\, a_{n} \otimes b_{n}\,\right>\,S^{\,-\, 1}_{F}p_{i} \otimes S^{\,-\, 1}_{G}q_{j}\\
&= \left(\sum\limits_{\,i \,=\, 1}^{\,\infty}\left<\,S^{\,-\, 1}_{F}p,\, p_{\,i}\,|\,a_{\,2},\, \cdots,\, a_{\,n}\,\right>_{1}S^{\,-\, 1}_{F}p_{\,i}\,\right)\otimes\\
&\hspace{1cm}\left(\sum\limits_{\,j \,=\, 1}^{\,\infty}\left<\,S^{\,-\, 1}_{G}q,\, q_{\,j}\,|\,b_{\,2},\, \cdots,\, b_{\,n}\,\right>_{2}\,S^{\,-\, 1}_{G}q_{\,j}\right)\\
&=\, S^{\,-\, 1}_{F}\,S_{F}\left(\,S^{\,-\, 1}_{F}\,p\,\right) \,\otimes\, S^{\,-\, 1}_{G}\,S_{G}\left(\,S^{\,-\, 1}_{G}\,q\,\right)\\
& \,=\, S^{\,-\, 1}_{F}\,p \,\otimes\, S^{\,-\, 1}_{G}\,q \,=\, S^{\,-\, 1}_{F \,\otimes\, G}\,(\,p \,\otimes\, q\,).
\end{align*}
Hence, the corresponding frame operator for \,$\Lambda$\, is \,$S^{\,-\, 1}_{F \,\otimes\, G}$. 
\end{proof}

\begin{theorem}
Let \,$\{\,p_{\,i}\,\}_{i \,=\,1}^{\infty}$\, be a frame associated to \,$\left(\,a_{\,2},\, \cdots,\, a_{\,n}\,\right)$\, for \,$H$\, with bounds \,$A,\,B$\, and \,$\{\,q_{\,j}\,\}_{j \,=\,1}^{\infty}$\, be a frame associated to \,$(\,b_{\,2},\, \cdots$, \,$b_{\,n}\,)$\, for \,$K$\, with bounds \,$C,\,D$\, having their corresponding frame operators \,$S_{F}$\, and \,$S_{G}$, respectively.\;Then \,$ \left\{\,\Delta_{i\,j} \,=\, \left(\,U_{1} \otimes U_{2}\,\right)\,\left(p_{\,i} \otimes q_{\,j}\right)\,\right\}^{\infty}_{i,\,j \,=\, 1}$\, is a frame associated to \,$(\,a_{\,2} \,\otimes\, b_{\,2},\, \,\cdots$,\, $a_{\,n} \,\otimes\, b_{\,n}\,)$\, for \,$H \,\otimes\, K$\, if and only if \,$U_{1} \,\otimes\, U_{2}$\, is an invertible operator on \,$H_{F} \,\otimes\, K_{G}$.
\end{theorem}

\begin{proof}
First we suppose that \,$U_{1} \,\otimes\, U_{2}$\, is an invertible on \,$H_{F} \,\otimes\, K_{G}$.\;Then by Theorem \ref{th1.1}, \,$U_{1}$\, and \,$U_{2}$\, are invertible on \,$H_{F}$\, and \,$K_{G}$, respectively.\;For each \,$p \,\in\, H_{F}$\, and \,$q \,\in\, K_{G}$, we have
\begin{equation}\label{eq1.05} 
\left\|\,p \,,\, a_{\,2},\, \cdots,\, a_{\,n}\,\right\|_{1} \,\leq\, \left\|\,U^{\,-\,1}_{1}\,\right\|\,\left\|\,U^{\,\ast}_{1}\,(\,p\,) \,,\, a_{\,2} \,,\, \cdots \,,\, a_{\,n}\,\right\|_{1},\;\text{and}
\end{equation}
\begin{equation}\label{eq1.5}
\left\|\,q \,,\, b_{\,2},\, \cdots,\, b_{\,n}\,\right\|_{2} \,\leq\, \left\|\,U^{\,-\, 1}_{2}\,\right\|\,\left\|\,U^{\,\ast}_{2}\,(\,q\,) \,,\, b_{\,2},\, \cdots,\, b_{\,n}\,\right\|_{2}.\hspace{.7cm}
\end{equation}
Now, for each  \,$p \,\otimes\, q \,\in\, H_{F} \,\otimes\, K_{G}$, we have
\begin{align}
&\sum\limits_{i,\, j \,=\, 1}^{\,\infty}\,\left|\,\left<\,p \,\otimes\, q \,,\, \left(\,U_{1} \,\otimes\, U_{2}\,\right)\,\left(\,p_{\,i} \,\otimes\, q_{\,j}\,\right) \,|\, a_{\,2} \,\otimes\, b_{\,2}, \,\cdots,\, a_{\,n} \,\otimes\, b_{\,n}\,\right>\,\right|^{\,2}\nonumber\\ 
&=\, \sum\limits_{i,\, j \,=\, 1}^{\,\infty}\,\left|\,\left<\,p \,\otimes\, q \,,\, U_{1}\,(\,p_{\,i}\,) \,\otimes\, U_{2}\,(\,q_{\,j}\,) \,|\, a_{\,2} \,\otimes\, b_{\,2}, \,\cdots,\, a_{\,n} \,\otimes\, b_{\,n}\,\right>\,\right|^{\,2}\nonumber\\
&= \sum\limits_{\,i \,=\, 1}^{\,\infty}\left|\,\left<\,U_{1}^{\,\ast}\,p,\, p_{\,i}\,|\,a_{\,2},\, \cdots,\, a_{\,n}\,\right>_{1}\,\right|^{\,2}\sum\limits_{\,j \,=\, 1}^{\,\infty}\left|\,\left<\,U_{2}^{\,\ast}\,q,\, q_{\,j} \,|\, b_{\,2},\, \cdots,\, b_{\,n}\,\right>_{2}\,\right|^{\,2}\label{eq1.6}\\
&\leq\, B\,D\,\left\|\,U_{1}^{\,\ast}\,(\,p\,) \,,\, a_{\,2},\, \cdots,\, a_{\,n}\,\right\|_{1}^{\,2}\,\left\|\,U_{2}^{\,\ast}\,(\,q\,) \,,\, b_{\,2},\, \cdots,\, b_{\,n}\,\right\|_{2}^{\,2}\nonumber\\
&[\;\text{since \,$\{\,p_{\,i}\,\}_{i \,=\,1}^{\infty}$\, is a frame associated to \,$\left(\,a_{\,2},\, \cdots,\, a_{\,n}\,\right)$, and}\nonumber\\
&\text{\,$\{\,q_{\,j}\,\}_{j \,=\,1}^{\infty}$\, is a frame associated to \,$\left(\,b_{\,2},\, \cdots,\, b_{\,n}\,\right)$\,}\;]\nonumber\\
&\leq\, B\,D\,\left\|\,U_{1}^{\,\ast}\,\right\|^{\,2}\,\left\|\,U_{2}^{\,\ast}\,\right\|^{\,2}\,\left\|\,p \,,\, a_{\,2},\, \cdots,\, a_{\,n}\,\right\|_{1}^{\,2}\,\left\|\,q \,,\, b_{\,2},\, \cdots,\, b_{\,n}\,\right\|_{2}^{\,2}\nonumber\\
&=\, B\,D\,\left\|\,U_{1} \,\otimes\, U_{2}\,\right\|^{\,2}\,\left\|\,p \,\otimes\, q \,,\, a_{\,2} \,\otimes\, b_{\,2} \,\cdots,\, a_{\,n} \,\otimes\, b_{\,n}\,\right\|^{\,2}.\nonumber 
\end{align}
On the other hand, from (\ref{eq1.6})
\begin{align*}
&\sum\limits_{i,\, j \,=\, 1}^{\,\infty}\left|\,\left<\,p \otimes q,\, \left(U_{1} \otimes U_{2}\right)\,\left(p_{\,i} \otimes q_{\,j}\right)\,|\,a_{\,2} \otimes b_{\,2}, \,\cdots,\, a_{\,n} \otimes b_{\,n}\,\right>\,\right|^{\,2}\\
&\geq\, A\,C\,\left\|\,U_{1}^{\,\ast}\,(\,p\,) \,,\, a_{\,2},\, \cdots,\, a_{\,n}\,\right\|_{1}^{\,2}\,\left\|\,U_{2}^{\,\ast}\,(\,q\,) \,,\, b_{\,2},\, \cdots,\, b_{\,n}\,\right\|_{2}^{\,2}\\
&\geq\, \dfrac{A\,C\left\|\,p,\, a_{\,2},\, \cdots,\, a_{\,n}\,\right\|_{1}^{\,2}\,\left\|\,q,\, b_{\,2},\, \cdots,\, b_{\,n}\,\right\|_{2}^{\,2}}{\left\|\,U^{\,-\,1}_{1}\,\right\|^{\,2}\,\left\|\,U^{\,-\,1}_{2}\,\right\|^{\,2}}\;[\;\text{by (\ref{eq1.05}) and (\ref{eq1.5})}\;]\\
&=\, \dfrac{A\,C}{\left\|\,\left(\,U_{1} \,\otimes\, U_{2}\,\right)^{\,-\, 1}\,\right\|^{\,2}}\,\left\|\,p \,\otimes\, q \,,\, a_{\,2} \,\otimes\, b_{\,2}, \,\cdots,\, a_{\,n} \,\otimes\, b_{\,n}\,\right\|^{\,2}.
\end{align*}
Therefore, \,$\left\{\,\Delta_{i\,j}\,\right\}_{i,\,j \,=\, 1}^{\,\infty}$\, is a frame associated to \,$\left(a_{\,2} \otimes b_{\,2},\, \,\cdots,\, a_{\,n} \otimes b_{\,n}\right)$\, for \,$H \,\otimes\, K$.\\

Conversely, suppose that \,$\left\{\,\Delta_{i\,j}\,\right\}_{i,\,j \,=\, 1}^{\,\infty}$\, is a frame associated to \,$(\,a_{\,2} \,\otimes\, b_{\,2},\, \,\cdots$, \,$a_{\,n} \,\otimes\, b_{\,n}\,)$\, for \,$H \,\otimes\, K$.\;Now, for each \,$p \,\otimes\, q \,\in\, H_{F} \,\otimes\, K_{G}$, 
\begin{align*}
&\sum\limits_{i,\, j \,=\, 1}^{\,\infty}\,\left<\,p \,\otimes\, q \,,\, \Delta_{i\,j} \,|\, a_{\,2} \,\otimes\, b_{\,2}, \,\cdots,\, a_{\,n} \,\otimes\, b_{\,n}\,\right>\,\Delta_{i\,j}\\
&=\, \sum\limits_{i,\, j \,=\, 1}^{\,\infty}\left<\,p \otimes q,\, U_{1}\,p_{\,i} \otimes U_{2}\,q_{\,j}\,|\,a_{\,2} \otimes b_{\,2}, \,\cdots,\, a_{\,n} \otimes b_{\,n}\,\right>\,\left(U_{1}\,p_{\,i} \otimes U_{2}\,q_{\,j}\right)\\
&= \sum\limits_{\,i \,=\, 1}^{\,\infty}\left<\,U_{1}^{\,\ast}p,\, p_{i}\,|\,a_{2},\, \cdots,\, a_{n}\,\right>_{1}U_{1}\,p_{i} \,\otimes\, \sum\limits_{\,j \,=\, 1}^{\,\infty}\left<\,U_{2}^{\,\ast}q,\, q_{j}\,|\,b_{2},\, \cdots,\, b_{n}\,\right>_{2}U_{2}\,q_{j}\\
&=\, U_{1}\,S_{F}\,U_{1}^{\,\ast}\,(\,p\,) \,\otimes\,  U_{2}\,S_{G}\,U_{2}^{\,\ast}\,(\,q\,)\\
& = \left(\,U_{1} \,\otimes\, U_{2}\,\right)\,\left(\,S_{F} \,\otimes\, S_{G}\,\right)\,\left(\,U^{\,\ast}_{1} \,\otimes\, U^{\,\ast}_{2}\,\right)\,(\,p \,\otimes\, q\,)\\
&=\, \left(\,U_{1} \,\otimes\, U_{2}\,\right)\,S_{F \,\otimes\, G}\,\left(\,U_{1} \,\otimes\, U_{2}\,\right)^{\,\ast}\,(\,p \,\otimes\, q\,).
\end{align*}
Hence, the frame operator for \,$\left\{\,\Delta_{i\,j}\,\right\}_{i,\,j \,=\, 1}^{\,\infty}$\, is \,$\left(U_{1} \,\otimes\, U_{2}\right)\,S_{F \,\otimes\, G}\,\left(U_{1} \,\otimes\, U_{2}\right)^{\,\ast}$\, and therefore it is invertible.\;Also, we know that \,$S_{F \,\otimes\, G}$\, is invertible and hence \,$U_{1} \,\otimes\, U_{2}$\, is invertible on \,$H_{F} \,\otimes\, K_{G}$.
\end{proof}

\section{Dual frame in tensor product of $n$-Hilbert spaces}

 In this section, dual frame in \,$n$-Hilbert spaces and their tensor product are discussed.

\begin{definition}\label{defi1} 
Let \,$\left\{\,p_{\,i}\,\right\}^{\,\infty}_{i \,=\, 1}$\, be a frame associated to \,$\left(\,a_{\,2} \,,\, \cdots \,,\, a_{\,n}\,\right)$\, for \,$H$. Then a frame \,$\left\{\,q_{\,i}\,\right\}^{\,\infty}_{i \,=\, 1}$\, associated to \,$\left(\,a_{\,2} \,,\, \cdots \,,\, a_{\,n}\,\right)$\, satisfying
\[p \,=\, \sum\limits^{\,\infty}_{i \,=\, 1}\, \left<\,p,\, q_{\,i} \,|\, a_{\,2},\, \cdots,\, a_{\,n}\,\right>_{1}\,p_{\,i}\; \;\;\forall\; p \,\in\, H\]
is called a dual frame or alternative dual frame associated to \,$\left(a_{\,2},\, \cdots,\, a_{\,n}\right)$\, of \,$\left\{\,p_{\,i}\,\right\}^{\,\infty}_{i \,=\, 1}$.
\end{definition}

\begin{theorem}\label{th3.1}
Let \,$\left\{\,p_{\,i}\,\right\}^{\,\infty}_{i \,=\, 1}$\, and \,$\left\{\,q_{\,i}\,\right\}^{\,\infty}_{i \,=\, 1}$\, be two Bessel sequences associated to \,$\left(\,a_{\,2} \,,\, \cdots \,,\, a_{\,n}\,\right)$\, in \,$H$.\;Then the following are equivalent:
\begin{itemize}
\item[$(i)$] $p \,=\, \sum\limits^{\,\infty}_{i \,=\, 1}\, \left<\,p,\, q_{\,i} \,|\, a_{\,2},\, \cdots,\, a_{\,n}\,\right>_{1}\,p_{\,i}; \;\forall\; p \,\in\, H_{F}$.
\item[$(ii)$] $p \,=\, \sum\limits^{\,\infty}_{i \,=\, 1}\, \left<\,p,\, p_{\,i} \,|\, a_{\,2},\, \cdots,\, a_{\,n}\,\right>_{1}\,q_{\,i}; \;\forall\; p \,\in\, H_{F}$.
\end{itemize}  
\end{theorem}

\begin{proof}$(\,i\,) \,\Rightarrow\, (\,ii\,)$
Let \,$T_{F}\; \;\text{and}\; \,T_{G}$\, be the pre-frame operators of \,$\left\{\,p_{\,i}\,\right\}^{\,\infty}_{i \,=\, 1}$\, and \,$\left\{\,q_{\,i}\,\right\}^{\,\infty}_{i \,=\, 1}$, respectively.\;Composing \,$T_{F}$\, with the adjoint of \,$T_{G}$, for all \,$p \,\in\, H_{F}$, we get 
\[T_{F}\,T_{G}^{\,\ast} \,:\, H_{F} \,\to\, H_{F},\; T_{F}\,T_{G}^{\,\ast}\,(\,p\,) \,=\, \sum\limits^{\,\infty}_{i \,=\, 1}\, \left<\,p,\, q_{\,i} \,|\, a_{\,2},\, \cdots,\, a_{\,n}\,\right>_{1}\,p_{\,i}.\]
Now, in terms of pre-frame operators \,$(\,i\,)$\, can be written as \;$T_{F}\,T_{G}^{\,\ast} \,=\, I_{F}$\; and this equivalent to \,$T_{G}\,T^{\,\ast}_{F} \,=\, I_{F}$\,.\;Therefore, for each \,$p \,\in\, H_{F}$, 
\[p \,=\, \sum\limits^{\,\infty}_{i \,=\, 1}\, \left<\,p,\, p_{\,i} \,|\, a_{\,2},\, \cdots,\, a_{\,n}\,\right>_{1}\,q_{\,i}.\]
Similarly, \,$(\,ii\,) \,\Rightarrow\, (\,i\,)$\, follows. 
\end{proof}

\begin{remark}
Suppose that the equivalent conditions of Theorem \ref{th3.1} are satisfied.\;Then using Cauchy-Schwartz inequality, for every \,$p \,\in\, H_{F}$, we have
\begin{align*}
&\left\|\,p,\, a_{\,2},\, \cdots,\, a_{\,n}\,\right\|_{1}^{\,2} \,=\, \left<\,p,\, p \,|\, a_{\,2},\, \cdots,\, a_{\,n}\,\right>_{1}\\
&\,=\, \left<\,\sum\limits^{\,\infty}_{i \,=\, 1}\, \left<\,p,\, p_{\,i} \,|\, a_{\,2},\, \cdots,\, a_{\,n}\,\right>_{1}\,q_{\,i},\, p \,|\, a_{\,2},\, \cdots,\, a_{\,n}\,\right>_{1}\\
&=\, \sum\limits^{\,\infty}_{i \,=\, 1}\, \left<\,p,\, p_{\,i} \;|\, a_{\,2},\, \cdots,\, a_{\,n}\,\right>_{1}\,\left<\,q_{\,i},\, p \,|\, a_{\,2},\, \cdots,\, a_{\,n}\,\right>_{1}\\
& \leq \left(\sum\limits^{\,\infty}_{i \,=\, 1}\left|\,\left<\,p,\, p_{\,i}\,|\,a_{2},\, \cdots,\, a_{n}\,\right>_{1}\,\right|^{\,2}\right)^{1 \,/\, 2} \left(\,\sum\limits^{\,\infty}_{i \,=\, 1}\left|\,\left<\,p,\, q_{\,i} \,|\, a_{2},\, \cdots,\, a_{n}\,\right>_{1}\,\right|^{\,2}\right)^{1 \,/\, 2}\\
&\leq\, \left(\,\sum\limits^{\,\infty}_{i \,=\, 1}\,\left|\,\left<\,p,\, p_{\,i} \,|\, a_{\,2},\, \cdots,\, a_{\,n}\,\right>_{1}\,\right|^{\,2}\,\right)^{1 \,/\, 2}\, B^{1 \,/\, 2}\; \left\|\,p,\, a_{\,2},\, \cdots,\, a_{\,n}\,\right\|_{1}\\
&[\;\text{since \,$\left\{\,q_{\,i}\,\right\}^{\,\infty}_{i \,=\, 1}$\, is a Bessel sequences associated to \,$\left(a_{\,2},\, \cdots,\, a_{\,n}\right)$\,}\;]\\
&\Rightarrow\, \dfrac{1}{B}\; \left\|\,p,\, a_{\,2},\, \cdots,\, a_{\,n}\,\right\|_{1}^{\,2} \,\leq\, \sum\limits^{\,\infty}_{i \,=\, 1}\,\left|\,\left<\,p,\, p_{\,i} \,|\, a_{\,2},\, \cdots,\, a_{\,n}\,\right>_{1}\,\right|^{\,2}.
\end{align*}
This shows that \,$\left\{\,p_{\,i}\,\right\}^{\,\infty}_{i \,=\, 1}$\, is a frame associated to \,$\left(\,a_{\,2},\, \cdots,\, a_{\,n}\,\right)$\, for \,$H$.\;Similarly, it can be shown that \,$\left\{\,q_{\,i}\,\right\}^{\,\infty}_{i \,=\, 1}$\, is also a frame associated to \,$\left(\,a_{\,2},\, \cdots,\, a_{\,n}\,\right)$\, for \,$H$.
\end{remark}

We now present the concept of a dual frame in \,$H \,\otimes\, K$.

\begin{remark}
Let \,$\left\{\,p_{\,i} \otimes q_{\,j}\,\right\}^{\,\infty}_{i,\,j \,=\, 1}$\, be a frame associated to \,$(\,a_{\,2} \otimes b_{\,2},\, \cdots$, \,$a_{\,n} \otimes b_{\,n}\,)$\, for \,$H \otimes K$.\,Then according to the definition \ref{defi1}, a frame \,$\left\{\,e_{\,i} \otimes h_{\,j}\,\right\}^{\,\infty}_{i,\,j \,=\, 1}$\, associated to \,$\left(\,a_{\,2} \,\otimes\, b_{\,2},\, \,\cdots,\, a_{\,n} \,\otimes\, b_{\,n}\,\right)$\, for \,$H \,\otimes\, K$\, satisfying 
\begin{equation}\label{eq2.1}
p \,\otimes\, q \,=\, \sum\limits_{i,\, j \,=\, 1}^{\,\infty}\,\left<\,p \,\otimes\, q,\, e_{\,i} \,\otimes\, h_{\,j} \,|\, a_{\,2} \,\otimes\, b_{\,2}, \,\cdots,\, a_{\,n} \,\otimes\, b_{\,n}\,\right>\,\left(\,p_{\,i} \,\otimes\, q_{\,j}\,\right),
\end{equation}
for all \,$p \,\otimes\, q \,\in\, H \,\otimes\, K$, can be consider as a dual frame associated to \,$(\,a_{\,2} \,\otimes\, b_{\,2},\, \cdots$, \,$a_{\,n} \,\otimes\, b_{\,n}\,)$\, of \,$\left\{\,p_{\,i} \,\otimes\, q_{\,j}\,\right\}^{\infty}_{i,\,j = 1}$.     
\end{remark}

\begin{remark}\label{note1}
According to the Theorem \ref{th2}, \,$\left\{\,p_{\,i} \,\otimes\, q_{\,j}\,\right\}^{\,\infty}_{i,\,j \,=\, 1}$\, and \,$\left\{\,e_{\,i} \,\otimes\, h_{\,j}\,\right\}^{\,\infty}_{i,\,j \,=\, 1}$\, are pair of dual frames associated to \,$(\,a_{\,2} \,\otimes\, b_{\,2},\, \,\cdots$, \,$a_{\,n} \,\otimes\, b_{\,n}\,)$\, for \,$H \,\otimes\, K$\, if and only if \,$\left\{\,p_{\,i} \,\otimes\, q_{\,j}\,\right\}^{\,\infty}_{i,\,j \,=\, 1}$\, and \,$\left\{\,e_{\,i} \,\otimes\, h_{\,j}\,\right\}^{\,\infty}_{i,\,j \,=\, 1}$\, are pair of dual frames for \,$H_{F} \,\otimes\, K_{G}$.  
\end{remark}

\begin{theorem}\label{th3.2}
Let \,$\{\,p_{\,i}\,\}_{i \,=\,1}^{\infty} \,,\, \left\{\,e_{\,i}\,\right\}^{\,\infty}_{i \,=\, 1}$\, be a pair of dual frames associated to \,$\left(\,a_{\,2},\, \cdots,\, a_{\,n}\,\right)$\, for \,$H$\, and \,$\{\,q_{\,j}\,\}_{j \,=\,1}^{\infty} \,,\, \left\{\,h_{\,j}\,\right\}^{\,\infty}_{j \,=\, 1}$\, be a pair of dual frames associated to \,$\left(\,b_{\,2},\, \cdots,\, b_{\,n}\,\right)$\, for \,$K$.\;Then \,$\left\{\,e_{\,i} \,\otimes\, h_{\,j}\,\right\}^{\,\infty}_{i,\,j \,=\, 1}$\, is a dual frame associated to \,$\left(\,a_{\,2} \,\otimes\, b_{\,2},\, \,\cdots,\, a_{\,n} \,\otimes\, b_{\,n}\,\right)$\, of \,$\left\{\,p_{\,i} \,\otimes\, q_{\,j}\,\right\}^{\,\infty}_{i,\,j \,=\, 1}$.     
\end{theorem}

\begin{proof}
By Theorem \ref{th2.1}, \,$\left\{\,p_{\,i} \,\otimes\, q_{\,j}\,\right\}^{\,\infty}_{i,\,j \,=\, 1}$, \,$\left\{\,e_{\,i} \,\otimes\, h_{\,j}\,\right\}^{\,\infty}_{i,\,j \,=\, 1}$\, are frames associated to \,$\left(\,a_{2} \,\otimes\, b_{2},\, \,\cdots,\, a_{n} \,\otimes\, b_{n}\,\right)$\, for \,$H \,\otimes\, K$.\;Since \,$\left\{\,e_{\,i}\,\right\}^{\,\infty}_{i \,=\, 1}$\, and \,$\left\{\,h_{\,j}\,\right\}^{\,\infty}_{j \,=\, 1}$\, are dual frames associated to \,$\left(a_{2},\, \cdots,\, a_{n}\right)$\, and \,$\left(b_{2},\, \cdots,\, b_{n}\right)$\, of \,$\{\,p_{\,i}\,\}_{i \,=\,1}^{\infty}$\, and \,$\{\,q_{\,j}\,\}_{j \,=\,1}^{\infty}$, respectively, for all \,$p \,\in\, H$, \,$q \,\in\, K$, 
\[p \,=\, \sum\limits^{\,\infty}_{i \,=\, 1}\left<\,p,\, e_{\,i} \,|\, a_{\,2},\, \cdots,\, a_{\,n}\,\right>_{1}\,p_{\,i},\, \;\text{and}\; q \,=\, \sum\limits^{\,\infty}_{j \,=\, 1}\left<\,q,\, h_{j} \,|\, b_{\,2},\, \cdots,\, b_{\,n}\,\right>_{2}\,q_{j}.\] 
Then, for all \,$p \,\otimes\, q \,\in\, H \,\otimes\, K$, we have
\begin{align*}
&p \,\otimes\, q \\
&\,=\, \left(\sum\limits^{\,\infty}_{i \,=\, 1} \left<\,p,\, e_{i}\,|\,a_{2},\, \cdots,\, a_{n}\,\right>_{1}p_{i}\right) \otimes \left( \sum\limits^{\,\infty}_{j \,=\, 1} \left<\,q,\, h_{j}\,|\,b_{2},\, \cdots,\, b_{n}\,\right>_{2}q_{j}\right)\\
&=\, \sum\limits_{i,\, j \,=\, 1}^{\,\infty}\,\left<\,p \,,\, e_{\,i} \;|\; a_{\,2} \,,\, \cdots \,,\, a_{\,n}\,\right>_{1}\, \left<\,q \,,\, h_{\,j} \;|\; b_{\,2} \,,\, \cdots \,,\, b_{\,n}\,\right>_{2}\,\left(\,p_{\,i} \,\otimes\, q_{\,j}\,\right)\\
&=\, \sum\limits_{i,\, j \,=\, 1}^{\,\infty}\,\left<\,p \,\otimes\, q,\, e_{\,i} \,\otimes\, h_{\,j} \,|\, a_{\,2} \,\otimes\, b_{\,2}, \,\cdots,\, a_{\,n} \,\otimes\, b_{\,n}\,\right>\,\left(\,p_{\,i} \,\otimes\, q_{\,j}\,\right).
\end{align*}
This completes the proof.    
\end{proof}

\begin{theorem}
Let \,$\{\,p_{\,i}\,\}_{i \,=\,1}^{\infty},\, \left\{\,e_{\,i}\,\right\}^{\,\infty}_{i \,=\, 1}$\, be a pair of dual frames associated to \,$\left(a_{\,2},\, \cdots,\, a_{\,n}\right)$\, for \,$H$\, and \,$\{\,q_{\,j}\,\}_{j \,=\,1}^{\infty},\, \left\{\,h_{\,j}\,\right\}^{\,\infty}_{j \,=\, 1}$\, be a pair of dual frames associated to \,$\left(b_{\,2},\, \cdots,\, b_{\,n}\right)$\, for \,$K$.\;Suppose \,$U \,\in\, \mathcal{B}\,(\,H_{F}\,)$\, and \,$V \,\in\, \mathcal{B}\,(\,K_{G}\,)$\, are unitary operators.\,Then \,$\Lambda \,=\, \left\{\left(U \otimes V\right)\left(p_{\,i} \otimes q_{\,j}\right)\right\}^{\,\infty}_{i,\,j = 1}$\, and \,$\Gamma \,=\, \left\{\,\left(\,U \,\otimes\, V\,\right)\,\left(\,e_{\,i} \,\otimes\, h_{\,j}\,\right)\,\right\}^{\,\infty}_{i,\,j \,=\, 1}$\, also form a pair of dual frames associated to \,$\left(\,a_{\,2} \,\otimes\, b_{\,2},\, \,\cdots,\, a_{\,n} \,\otimes\, b_{\,n}\,\right)$\, for \,$H \,\otimes\, K$.   
\end{theorem}

\begin{proof}
By Theorem \ref{th3.2}, \,$\left\{\,p_{\,i} \,\otimes\, q_{\,j}\,\right\}^{\,\infty}_{i,\,j \,=\, 1}$\, and \,$\left\{\,e_{\,i} \,\otimes\, h_{\,j}\,\right\}^{\,\infty}_{i,\,j \,=\, 1}$\, form a pair of dual frames associated to \,$\left(\,a_{\,2} \,\otimes\, b_{\,2},\, \,\cdots,\, a_{\,n} \,\otimes\, b_{\,n}\,\right)$\, for \,$H \,\otimes\, K$.\;Now, for each \,$p \,\otimes\, q \,\in\, H_{F} \,\otimes\, K_{G}$, we have
\begin{align*}
&\sum\limits_{i,\, j \,=\, 1}^{\,\infty}\,\left|\,\left<\,p \,\otimes\, q,\, \left(\,U \,\otimes\, V\,\right)\,\left(\,e_{\,i} \,\otimes\, h_{\,j}\,\right) \,|\, a_{\,2} \,\otimes\, b_{\,2}, \,\cdots,\, a_{\,n} \,\otimes\, b_{\,n}\,\right>\,\right|^{\,2}\\
&=\, \sum\limits_{i,\, j \,=\, 1}^{\,\infty}\,\left|\,\left<\,p \,\otimes\, q,\, \left(\,U\,e_{\,i} \,\otimes\, V\,h_{\,j}\,\right) \,|\, a_{\,2} \,\otimes\, b_{\,2}, \,\cdots,\, a_{\,n} \,\otimes\, b_{\,n}\,\right>\,\right|^{\,2}\\
&= \sum\limits_{\,i \,=\, 1}^{\,\infty}\left|\,\left<\,U^{\,\ast}p,\, e_{i}\,|\,a_{2},\, \cdots,\, a_{n}\,\right>_{1}\right|^{\,2} \otimes \sum\limits_{\,j \,=\, 1}^{\,\infty}\left|\,\left<\,V^{\,\ast}q,\, h_{j}\,|\,b_{2},\, \cdots,\, b_{n}\,\right>_{2}q_{j}\,\right|^{\,2}.
\end{align*}
Since \,$\{\,e_{\,i}\,\}_{i \,=\,1}^{\infty}$\, is a frame associated to \,$\left(\,a_{\,2},\, \cdots,\, a_{\,n}\,\right)$\, for \,$H$\, and \,$\{\,h_{\,j}\,\}_{j \,=\,1}^{\infty}$\, is a frame associated to \,$\left(\,b_{\,2},\, \cdots,\, b_{\,n}\,\right)$\, for \,$K$, the above calculation shows that \,$\Gamma$\, is a frame associated to \,$\left(\,a_{\,2} \,\otimes\, b_{\,2},\, \,\cdots,\, a_{\,n} \,\otimes\, b_{\,n}\,\right)$\, for \,$H \,\otimes\, K$.\;Similarly, it can be shown that \,$\Lambda$\, is a frame associated to \,$\left(\,a_{\,2} \,\otimes\, b_{\,2},\, \,\cdots,\, a_{\,n} \,\otimes\, b_{\,n}\,\right)$\, for \,$H \,\otimes\, K$.  
Furthermore, for each \,$p \,\otimes\, q \,\in\, H_{F} \,\otimes\, K_{G}$, we have
\begin{align*}
&\sum\limits_{i,\, j = 1}^{\infty}\left<\,p \otimes q,\, \left(U \otimes V\right)\,\left(e_{i} \otimes h_{j}\right)\,|\,a_{2} \otimes b_{2}, \,\cdots,\, a_{n} \otimes b_{n}\,\right>\,\left(U \otimes V\right)\left(p_{i} \otimes q_{j}\right)\\
&=\, \sum\limits_{i,\, j \,=\, 1}^{\,\infty}\left<\,p \,\otimes\, q,\, \left(\,U\,e_{i} \,\otimes\, V\,h_{j}\,\right) \,|\, a_{2} \,\otimes\, b_{2}, \,\cdots,\, a_{n} \,\otimes\, b_{n}\,\right>\,\left(\,U\,p_{i} \,\otimes\, V\,q_{j}\,\right)\\
&= U\sum\limits_{\,i = 1}^{\infty}\left<\,U^{\,\ast}p,\, e_{\,i}\,|\,a_{\,2},\, \cdots,\, a_{\,n}\,\right>_{1}p_{\,i} \,\otimes\, V\sum\limits_{\,j = 1}^{\,\infty}\left<\,V^{\,\ast}q,\, h_{\,j}\,|\,b_{\,2},\, \cdots,\, b_{\,n}\,\right>_{2}q_{\,j}\\
& \,=\, U\,U^{\,\ast}\,(\,p\,) \,\otimes\, V\,V^{\,\ast}\,(\,q\,) \,=\, p \,\otimes\, q\; \;[\;\text{since \,$U,\, V$\, are unitary operators}\;].\\
&[\;\text{Since \,$\{\,p_{i}\,\}_{i \,=\,1}^{\infty}$, \,$\left\{\,e_{i}\,\right\}^{\,\infty}_{i \,=\, 1}$\, are dual frames associated to \,$(\,a_{2},\, \cdots,\, a_{n}\,)$,}\\
&\text{ and\,$\{\,q_{\,j}\,\}_{j \,=\,1}^{\infty},\, \left\{\,h_{\,j}\,\right\}^{\,\infty}_{j \,=\, 1}$\, are dual frames associated to \,$\left(\,b_{\,2},\, \cdots,\, b_{\,n}\,\right)$}\;].
\end{align*}
Hence, according to the remark \ref{note1}, \,$\Lambda$\, and \,$\Gamma$\, form a pair of dual frames associated to \,$\left(\,a_{\,2} \,\otimes\, b_{\,2},\, \,\cdots,\, a_{\,n} \,\otimes\, b_{\,n}\,\right)$\, for \,$H \,\otimes\, K$.
\end{proof}


\bigskip
\noindent
{\bf Acknowledgment.}
The authors would like to thank the editor and the referees for their helpful suggestions and comments to improve this paper.

\bibliographystyle{amsplain}

\vspace{0.1in}
\hrule width \hsize \kern 1mm
\end{document}